\documentclass{amsart} 
\usepackage{amssymb}
\usepackage{amsmath}
\usepackage{euscript}
\begin{document} 
\input epsf.sty

\title{Prime and composite Laurent polynomials}  
\author{F. Pakovich}
\address{Department of Mathematics,
Ben Gurion University 
P.O.B. 653, Beer Sheva 84105, Israel}
\email{pakovich@math.bgu.ac.il}
\date{}

\subjclass{Primary 30D05; Secondary 14H30}

\keywords{Ritt's theorems, decompositions of rational functions, decompositions of Laurent polynomials}

\begin{abstract}
In the paper \cite{r1} Ritt constructed the theory of functional
decompositions of polynomials with complex coefficients. In particular, he  
described explicitly polynomial solutions of  
the functional equation $f(p(z))=g(q(z)).$ In this paper we 
study the equation above in the case when $f,g,p,q$ are holomorphic functions on compact Riemann surfaces. We also 
construct a self-contained theory of functional decompositions of rational functions 
with at most two poles generalizing the Ritt theory. In particular, we 
give new proofs of the theorems of 
Ritt and of the theorem of Bilu and Tichy.

\end{abstract}

\maketitle

\def\be{\begin{equation}}
\def\ee{\end{equation}}
\def\bs{$\square$ \vskip 0.2cm}
\def\d{{\rm d}} 
\def\D{{\rm D}} 
\def\I{{\rm I}} 
\def\C{{\mathbb C}} 
\def\N{{\mathbb N}} 
\def\P{{\mathbb P}}
\def\Z{{\mathbb Z}}
\def\R{{\mathbb R}} 
\def\ord{{\rm ord}}
\def\f{\EuScript}

\def\e{\eqref}
\def\phi{{\varphi}}
\def\v{{\varepsilon}} 
\def\deg{{\rm deg\,}}
\def\Aut{{\rm Aut\,}}
\def\Det{{\rm Det}}
\def\dim{{\rm dim\,}} 
\def\Ker{{\rm Ker\,}} 
\def\Gal{{\rm Gal\,}}
\def\St{{\rm St\,}} 
\def\exp{{\rm exp\,}} 
\def\cos{{\rm cos\,}} 
\def\diag{{\rm diag\,}} 
\def\GCD{{\rm GCD }}
\def\LCM{{\rm LCM }}
\def\mod{{\rm mod\ }}
\def\mm{{\rm mult}}
\def\c{\circ}

\def\bp{\begin{proposition}}
\def\ep{\end{proposition}}
\def\bt{\begin{theorem}}
\def\et{\end{theorem}}
\def\be{\begin{equation}}
\def\bee{\begin{equation*}}
\def\la{\label}
\def\l{\lambda}
\def\ee{\end{equation}}
\def\eee{\end{equation*}}
\def\bl{\begin{lemma}}
\def\el{\end{lemma}}
\def\bc{\begin{corollary}}
\def\ec{\end{corollary}}
\def\pr{\noindent{\it Proof. }}
\def\note{\noindent{\bf Note. }}
\def\bd{\begin{definition}}
\def\ed{\end{definition}}
\def\qed{$\ \ \Box$ \vskip 0.2cm}

\newtheorem{theorem}{Theorem}[section]
\newtheorem{lemma}[theorem]{Lemma}
\newtheorem{corollary}[theorem]{Corollary}
\newtheorem{proposition}[theorem]{Proposition}

\section{Introduction} 
Let $F$ be a rational function with complex coefficients.
The function $F$ is called {\it indecomposable} if
the equality $F=F_2\circ F_1$, where $F_2\circ F_1$ denotes a superposition $F_2(F_1(z))$ of rational functions $F_1,$ $F_2,$ implies that at least one of the functions $F_1,F_2$ is of degree one. 
Any representation 
of a rational function $F$ in the form $F=F_r\circ F_{r-1}\circ \dots \circ F_1,$
where $F_1,F_2,\dots, F_r$ are rational functions,
is called {\it a decomposition} of $F.$ A decomposition is called {\it maximal}
if all $F_1,F_2,\dots, F_r$ are
indecomposable and of degree greater than one.

In general, a rational function may have many maximal decompositions 
and the ultimate goal
of the decomposition theory of rational functions is to
describe the general structure of all maximal
decompositions up to an equivalence, where by definition
two decompositions having an equal number of terms
$$F=F_r\circ F_{r-1}\circ \dots \circ F_1 \ \ \ \ {\rm and} \ \ \ \  
F=G_{r}\circ G_{r-1}\circ \dots \circ G_1$$
are called equivalent if either $r=1$ and $F_1=G_1$, or $r\geq 2$ and there exist rational functions $\mu_i,$ $1\leq i \leq r-1,$ of degree 1 such that 
$$F_r=G_r\circ \mu_{r-1}, \ \ \ 
F_i=\mu_{i}^{-1}\circ G_i \circ \mu_{i-1}, \ \ \ 1<i< r, \ \ \ {\rm and} \ \ \ F_1=\mu_{1}^{-1}\circ G_1.
$$
Essentially, the unique class of rational functions for which this problem is completely solved 
is the class of polynomials investigated by Ritt in his classical paper \cite{r1}.

The results of Ritt can be summarized in the form
of two theorems usually called the first and the second Ritt
theorems (see \cite{r1}, \cite{sch}). 
The first Ritt theorem states that any two maximal decompositions $\f D, \f E$ 
of a 
polynomial $P$ have an equal number of terms and there exists 
a chain of maximal decompositions $\f F_i$, $1\leq i \leq s,$ of $P$ such that 
$\f F_1=\f D,$ $\f F_s\sim \f E,$ and $\f F_{i+1}$ is obtained from $\f F_i$ 
by replacing 
two successive functions $A\circ C$ in $\f F_i$
by two other functions $B\circ D$ such that 
\be \la{-0} A\circ C=B\circ D.\ee

The second Ritt theorem states that if $A,B,C,D$ is a polynomial solution 
of \eqref{-0} such that 
$${\rm GCD}(\deg A, \deg B)=1, \ \ \ {\rm GCD}(\deg C,\deg D)=1$$ (this condition is satisfied in particular if $A,B,C,D$ are indecomposable)
then there exist polynomials $\tilde A,$ $\tilde B,$ $\tilde C,$ $\tilde D,$
$\mu_1,$ $\mu_2,$ where $\deg \mu_1=1,$ $\deg \mu_2=1,$
such that
$$A=\mu_1 \circ \tilde A, \ \ \ B=\mu_1 \circ \tilde B,
\ \ \ C=\tilde C \circ \mu_2, \ \ \ D=\tilde D \circ \mu_2$$ and either
$$\tilde A\circ\tilde C \sim T_n\circ T_m, \ \ \ \ \  
\tilde B\circ \tilde D \sim T_m\circ T_n,$$ where $T_m, T_n$ are the corresponding Chebyshev polynomials with $n,m\geq 1$ and $\GCD(n,m)=1,$ or
$$\tilde A\circ \tilde C \sim z^n\circ z^rR(z^n), \ \ \ \ \ 
\tilde  B\circ \tilde D \sim z^rR^n(z)
\circ z^n,$$ where $R$ is a polynomial, $r\geq 0,$ $n\geq 1,$ and $\GCD(n,r)=1$. Actually, 
the second 
Ritt theorem essentially remains true for arbitrary polynomial solutions of \eqref{-0}. 
The only difference in the formulation is that for the degrees of polynomials $\mu_1,$ $\mu_2$ in this case the equalities $$\deg \mu_1=\GCD(\deg A,\deg B),\ \ \ \deg \mu_2=\GCD(\deg C,\deg D)$$ hold
(see \cite{en}, \cite{tor}). 
Notice that an analogue of the second Ritt theorem 
holds also when the ground field is distinct from $\C$ (see \cite{za}).

For arbitrary rational functions the first Ritt
theorem fails to be true. Furthermore, there exist rational functions
having maximal decompositions of different length. The simplest 
examples of such functions 
can be constructed with the use of   
rational functions which are Galois coverings. These functions,
for the first time calculated by Klein in his
famous book \cite{klein}, are related to the finite subgroups $C_n,$ $D_n,$ $A_4,$ $S_4,$ $A_5$ of 
$\Aut\C\P^1$ and 
nowadays can be interpreted 
as Belyi functions of Platonic solids (see \cite{cg}, \cite{zv}).
Since for such a function $f$ 
its maximal decompositions  
correspond to maximal chains of subgroups 
of its monodromy group $G$, in order to find 
maximal decompositions of different length of $f$ 
it is enough to find the corresponding chains of subgroups of $G$, and it is not hard to check that 
for the groups $A_4,$ $S_4,$ and $A_5$ such chains exist (see \cite{gs}, \cite{mp2}).

The analogues of the second Ritt theorem for arbitrary rational solutions of equation \eqref{-0}
are known only in several cases. 
Let us mention some of them.
First, notice that the description of rational solution of \eqref{-0}
under condition that $C$ and $D$ are polynomials turns out to be quite simple
and substantially reduces to the description
of polynomial solutions of \eqref{-0} (see \cite{p2}). On the other hand, the problem of
description of rational solutions of \eqref{-0} under condition that $A$ and $B$ are polynomials is equivalent to the problem of description of algebraic curves of the form \be \la{cur} A(x)-B(y)=0,\ee  having a factor of genus zero, together with corresponding parametrizations.
A complete list of such curves is known only in the case when the corresponding factor has at most two points at infinity. In this case the problem
is closely related to the number theory and was studied first in the paper of Fried \cite{f1} and then in the papers of Bilu \cite{bilu1} and Bilu and Tichy \cite{bilu}. In particular, in \cite{bilu}
an explicit list of such curves, defined over any field of characteristic zero, was obtained.
Notice that the results of
\cite{f1}, \cite{bilu} generalize the second Ritt theorem since polynomial solutions of \eqref{-0} correspond to curves \eqref{cur}
having a factor of genus zero with one point at infinity.
Rational solutions of the equation
\be \la{az} A\circ C=A\circ D,\ee under condition that $A$ is a polynomial were described in \cite{az} (notice also the paper \cite{r2} where some partial results about equation \eqref{az} under condition that $A$ is a
rational function were obtained). Finally, a description of permutable rational functions was obtained in \cite{r3}
(see also \cite{er}).
Note that beside of connections with the number theory equation \eqref{-0} has also important connections with different 
branches of analysis 
(see e.g. recent papers \cite{p}, \cite{p2}, \cite{p1}, \cite{pppp}, \cite{pr1}).

In this paper we study the equation \be \la{df} h=f\circ p=g\circ q,\ee where $f:\, C_1\rightarrow \C\P^1,$ $g:\, C_2\rightarrow \C\P^1$ are fixed 
holomorphic
functions on fixed connected compact Riemann surfaces $C_1,C_2$ and $h:\, C\rightarrow \C\P^1,$
$p:\, C\rightarrow C_1,$ $q:\, C\rightarrow C_2$ are 
unknown holomorphic functions on unknown connected compact Riemann surface $C.$ 
We also apply the results obtained to equation \eqref{-0} 
with rational $A,B,C,D$ and on this base  
construct a self-contained decomposition theory of rational functions with at most {\it two} poles generalizing the Ritt theory. In particular, we prove analogues of Ritt theorems for such functions and 
reprove in a uniform way previous related results of \cite{r1}, \cite{f1}, \cite{bilu1}, \cite{bilu}.

Let $S\subset \C\P^1$ be a finite set and $z_0\in \C\P^1\setminus S.$ 
Our approach to equation \eqref{df} is based on the correspondence between
pairs consisting of a covering $f$ of $\C\P^1$, non-ramified outside of $S$, together with a 
point from $f^{-1}\{z_0\}$  
and subgroups of finite index in $\pi_1(\C\P^1\setminus S,z_0).$
The main advantage of the consideration of such pairs and 
subgroups, rather than just of functions and their monodromy groups,
is due to the fact that 
for any 
subgroups of finite index $A,B$ in $\pi_1(\C\P^1\setminus S,z_0)$ 
the subgroups $A\cap B$ and $<A,B>$ 
also are subgroups of finite index in $\pi_1(\C\P^1\setminus S,z_0)$ 
and we may transfer these operations to the corresponding pairs. 
The detailed description of the content of the paper is given below.

In Section \ref{s1} we describe the general structure of solutions of equation \eqref{df}.
We show (Theorem \ref{p2}) that 
there exists a finite number $o(f,g)$ of solutions $h_j,p_j,q_j$ of \eqref{df} 
such that any other solution may be obtained from them
and describe explicitly the monodromy of $h_j$ via the monodromy of $f,g.$ 
Furthermore, we show (Proposition \ref{compon}) that if $f,$ $g$ are rational functions then the Riemann surfaces on which 
the functions $h_j,$ $1\leq j \leq o(f,g),$ are defined may be identified with irreducible components of the algebraic curve $f(x)-g(y)=0$. In particular, being applied to 
polynomials $A,B$ our construction provides a criterion 
for the irreducibi\-lity of
curve \eqref{cur} via the monodromy groups of $A$ and $B$ 
useful for applications (see e.g. \cite{pppp}).

By the analogy with rational functions we will call a pair of holomorphic functions $f,g$
irreducible if $o(f,g)=1$. In Section \ref{s2} we study properties of irreducible and reducible pairs. In particular, we give  
a criterion (Theorem \ref{p5}) for a pair $f,g$ to be irreducible in terms of the corresponding subgroups of 
$\pi_1(\C\P^1\setminus S,z_0)$ and establish 
the following result about reducible pairs generalizing the corresponding result of Fried \cite{f2} about rational functions (Theorem \ref{p6}): if a pair of holomorphic functions $f,g$ is reducible 
then there exist holomorphic functions $\tilde f,$ $\tilde g,$ $p,$ $q$ 
such that
$$ f= \tilde f\circ p, \ \ \ g=\tilde g\circ q, \ \ \ o(f,g)=o(\tilde f, \tilde g), $$ 
and the Galois closures of $\tilde f$ and $\tilde g$ coincide. We also show (Theorem \ref{rit2}) that if in \eqref{df} the pair $f,$ $g$ is irreducible then the indecomposability of $q$ implies 
the indecomposability of $f$. Notice that the last result turns out to be quite useful for 
applications related to possible generalizations of the first Ritt theorem
(see Section \ref{rita}).

Further, in Section \ref{s3} we 
study properties of equation \eqref{df} in the case when
$f,g$
are ``generalized polynomials'' that is holomorphic functions for which the preimage of infinity 
contains a unique point.
In particular, we establish the following, highly useful for the study of equation \eqref{-0},
result (Corollary \ref{copo}): 
if $A,B$ are polynomials of the same degree and $C,D$ are rational functions such that equality \eqref{-0} holds then there exist a
rational function $W$, mutually distinct points of the complex sphere $\gamma_i,$ $1\leq i \leq r,$ and complex numbers 
$\alpha_i, \beta_i$ $0\leq i \leq r,$
such that 
$$C=\left(\alpha_0+\frac{\alpha_1}{z-\gamma_1}+\dots + \frac{\alpha_r}{z-\gamma_r}\right)
\circ W, \ \ \ \ \ D= \left(\beta_0+\frac{\beta_1}{z-\gamma_1}+\dots + \frac{\beta_r}{z-\gamma_r}\right)\circ W.$$

In Section \ref{rita} we propose an approach to possible generalizations of the first Ritt theorem
to more wide than polynomials classes of functions. We
introduce the conception of a closed class of rational functions as of a subset $\f R$ of $\C(z)$
such that the condition $G\circ H \in \f R$ implies that $G\in \f R$, $H\in \f R.$ 
The prototypes for this definition are closed 
classes $\f R_k,$ $k\geq 1,$ consisting of rational functions $F$ 
for which \be \la{mini} \min_{z\in \C\P^1}\vert  F^{-1}\{z\}\vert \leq k, \ee 
where $\vert F^{-1}\{z\} \vert$ denotes the cardinality of the set $F^{-1}\{z\}$.
Notice that 
since for any $F\in \f R_1$ there exist rational functions $\mu_1,$ $\mu_2$ of degree 1 such that $\mu_1\circ F\circ \mu_2$ is a 
polynomial,
the Ritt theorems can be interpreted as a decomposition theory for the class $\f R_1.$  
The main result of Section \ref{rita} (Theorem \ref{rit1}) states that in order to check that the first Ritt theorem holds for maximal decompositions 
of rational functions from a closed class $\f R$
it is enough to check that it holds for a certain subset of maximal decompositions 
which is considerably smaller than the whole set.
For example, for the class $\f R_1$ this subset turns out to be empty that provides  
a new proof of the first Ritt for this class (Corollary \ref{rit3}). Later, in Section 9, using 
this method we also 
show that the first Ritt theorem 
remains true for 
the class $\f R_2$.

In the rest of the paper, using the results obtained,
we construct explicitly the decomposition theory for the class $\f R_2$.
There are several reasons which
make the problem interesting. First, since $\f R_1\subset \f R_2$,
the decomposition theory for $\f R_2$ is a natural 
generalization of the Ritt theory. Furthermore, 
the equation 
\be \la{ld} L=A\circ C=B\circ D,\ee where $L\in \f R_2$ and 
$A,B,C,D$ are rational functions, is closely related to 
the equation \be \la{ppoo} h=A\circ f=B\circ g,\ee where $A,B$ are rational functions while $h,f,g$ are 
entire transcendental functions and the description of solutions of \eqref{ld}
yields a description of solutions of \eqref{ppoo} (see \cite{pr1}). Finally, notice that 
polynomials solutions of \eqref{-0}
naturally appear in the study of the polynomial moment problem which arose recently in connection with the ``model'' problem for the 
Poincare center-focus problem (see e. g. \cite{p}, \cite{y}). The corresponding moment problem for Laurent polynomials, which is related to the Poincare problem even to a greater extent than the polynomial moment problem, 
is still open and the decomposition theory for $\f R_2$ can be considered as a preliminary
step in the investigation of this problem. 

It was observed by the author several years ago 
that the description of ``double decompositions'' \eqref{ld} of functions from $\f R_2$ 
(``the second Ritt theorem'' for $\f R_2$) mostly reduces to the classification of curves \eqref{cur} having a factor of genus
0 with at most two points at infinity. Indeed, without loss of generality we may assume that 
the minimum in \eqref{mini} attains at infinity and that $L^{-1}\{\infty\}\subseteq \{0,\infty\}$.
In other words, we may assume that $L$ is a Laurent polynomial. Further, it follows easily from the condition  $L^{-1}\{\infty\}\subseteq \{0,\infty\}$ that any decomposition $U\circ V$ of $L$ 
is equivalent either to a decomposition $A\circ L_1$, where $A$ is a polynomial and $L_1$ is a Laurent polynomial, or to a decomposition $L_2\circ B$, where $L_2$ is
a Laurent polynomial and $B=cz^d$ for some $c\in \C$ and $d\geq 1.$
Therefore, the 
description of double decompositions of functions
from $\f R_2$
reduces to the solution of the following three equations:
\be \la{1} A\circ L_1=B\circ L_2
\ee where 
$A, B$ are polynomials and $L_1,L_2$ are Laurent polynomials, \be \la{2} A\circ L_1=L_2\circ z^d,\ee where 
$A$ is a polynomial and $L_1,L_2$ are Laurent polynomials, and
\be \la{3}L_1\circ z^{d_1}= L_2\circ z^{d_2},
\ee where $L_1,L_2$ are Laurent polynomials.
Observe now that if $A,B,L_1,L_2$ is a solution of equation \eqref{1} then 
corresponding curve \eqref{cur} has a factor
of genus $0$ with at most two points at infinity and vice versa for any such a curve the corresponding factor may be parametrized by some Laurent polynomials providing a solution of \eqref{1}.
Therefore, the description of solutions of equation \eqref{1} essentially reduces to the description of curves \eqref{cur} having a factor
of genus $0$ with at most two points at infinity.
On the other hand, equations \eqref{2} and \eqref{3} turn out to be much easier for the analysis in view of the presence of symmetries.

Although the 
result of 
Bilu and Tichy obtained in the paper 
\cite{bilu} (which in its turn uses the results of 
the papers \cite{bilu1}, \cite{f1}, \cite{f2}) reduces the solution of equation \eqref{1}
to an elementary problem of finding of parameterizations of the corresponding curves, 
in this paper we give an independent analysis of this equation
in view of the following reasons.
First, we wanted to provide a self contained 
exposition of the decomposition theory for the class $\f R_2$
since we believe that such an exposition may be interesting  
for the wide audience. 
Second,  
our approach contains some new ideas and by-product results which seem to be interesting by themselves.

Our analysis of equations \eqref{1}, \eqref{2}, \eqref{3} splits into three parts. 
In Section \ref{urrr} using Corollary \ref{copo} we solve equations \eqref{2}, \eqref{3}.
In Section \ref{reduc} using Theorem \ref{p6} combined with Corollary \ref{copo}
we show (Theorem \ref{red}) 
that equation \eqref{1} in the case when curve \eqref{cur} is reducible reduces either to the irreducible case 
or to the case when $$A\circ L_1=B\circ L_2=\frac{1}{2}\left(z^d+\frac{1}{z^d}\right), \ \ \ d>1.$$ 
Finally, in Section \ref{irreduc} we solve equation \eqref{1} in the case when
curve \eqref{cur} is irreducible. Our approach to this case 
is similar to the one used in the paper \cite{bilu} and consists of the analysis of 
the condition that the genus $g$ of \eqref{cur} is zero. However, we use a different form of the formula for $g$ and replace the 
conception of ``extra'' points which goes back to Ritt by a more transparent conception. 

Eventually, in Section \ref{last} of the paper, as a corollary of the classification of double decompositions of functions from $\f R_2$ 
and Theorem \ref{rit1},
we show (Theorem \ref{rittt}) that the first Ritt theorem extends
to the class $\f R_2$. 
The results of the paper concerning decompositions of functions from $\f R_2$ can be summarized in the form of the following theorem which includes in particular the Ritt theorems and the classifications of curves \eqref{cur} having a factor of genus 0 with two points at infinity.

\bt \la{1.1} Let 
$$L=A\circ C=B\circ D$$ be two decompositions of a rational function 
$L\in \f R_2$ into compositions of rational functions $A,C$ and $B,D$. 
Then 
there exist rational functions $R,W,\tilde A,\tilde B,\tilde C,\tilde D \in \f R_2$
such that   
$$ A=R \circ \tilde A, \ \ \ \ B=R \circ \tilde B, \ \ \ \ 
C=\tilde C \circ W, \ \ \ \ D=\tilde D \circ W,\ \ \ \ \tilde A\circ \tilde C=\tilde B\circ \tilde D$$
and, up to a possible replacement of $A$ by $B$ and $C$ by $D$, one of the following conditions holds:
\vskip 0.01cm 
$$\tilde A\circ \tilde C\sim z^n \circ z^rL(z^n),  \ \ \ \ \ \ 
\tilde B\circ \tilde D\sim  z^rL^n(z) \circ z^n,\leqno 1) $$ 
where $L$ is a Laurent polynomial, $r\geq 0,$ $n\geq 1,$  and 
$\GCD(n,r)=1;$ 
\vskip 0.01cm
$$\tilde A\circ \tilde C\sim  z^2 \circ \frac{z^2-1}{z^2+1}
S\left(\frac{2z}{z^2+1}\right),  \ \ \ \ \ \ \tilde B\circ \tilde D\sim  (1-z^2)S^2(z)\circ \frac{2z}{z^2+1},\leqno 2)$$
where $S$ is a polynomial;
\vskip 0.01cm 
$$\tilde A\circ \tilde C\sim T_n \circ T_m, \ \ \ \ \ \ \tilde B\circ \tilde D\sim T_m \circ T_n,\leqno 3)$$ 
where $T_n,T_m$ are the corresponding Chebyshev polynomials with $m,n\geq 1,$ and $\GCD(n,m)=1;$ 
\vskip 0.01cm
$$\tilde A\circ \tilde C\sim T_{n} \circ \frac{1}{2}\left(z^{m}+
\frac{1}{z^{m}}\right), \ \ \ \ \ \ 
\tilde B\circ \tilde D\sim  \frac{1}{2}\left(z^{m}
+\frac{1}{z^{m}}\right) \circ z^{n},\leqno 4)$$ 
where $m,n\geq 1$ and $\GCD(n,m)=1$;
\vskip 0.01cm
$$\tilde A\circ \tilde C\sim -T_{nl} \circ \frac{1}{2}\left(\varepsilon z^{m}+
\frac{1}{\varepsilon z^{m}}\right), \ \ \ \ \ \ 
\tilde B\circ \tilde D\sim T_{ml} \circ  \frac{1}{2}\left(z^{n}
+\frac{1}{z^{n}}\right),\leqno 5)$$    
where $T_{nl},T_{ml}$ are the corresponding Chebyshev polynomials with $m,n\geq 1,$ $l>1$, $\varepsilon^{nl}=-1$,
and $\GCD(n,m)=1$; 
\vskip -0.2cm 
$$\tilde A\circ \tilde C\sim (z^2-1)^3\circ 
\frac{3(3z^4+4z^3-6z^2+4z-1)}{(3z^2-1)^2}, \leqno 6) $$
$$\ \ \ \ \ \ \tilde B\circ \tilde D\sim (3z^4-4z^3)\circ \frac{4(9z^6-9z^4+18z^3-15z^2+6z-1)}{(3z^2-1)^3}.$$
\vskip 0.05cm
Furthermore, if $\f D,\f E$ are two maximal decompositions 
of $L$ then there exists 
a chain of maximal decompositions $\f F_i$, $1\leq i \leq s,$ of $L$ such that 
$\f F_1=\f D,$ $\f F_s\sim \f E,$ and $\f F_{i+1}$ is obtained from $\f F_i$ by replacing 
two successive functions in $\f F_i$ by two other functions with the same composition.
\et

\vskip 0.2cm
\noindent{\bf Acknowledgments}. The results of this paper were obtained mostly during the visits of the author to the Max-Planck-Institut f\"ur Mathematik in Summer 2005 and Spring 2007 and were partially announced in \cite{ar}. Seizing an opportunity the author would 
like to thank the Max-Planck-Institut for the hospitality. Besides, the author is grateful to Y. Bilu, M. Muzychuk, and M. Zieve for discussions of ideas and results of
this paper before its publication.

\section{\la{s1} Functional equation $h=f\circ p=g\circ q$}
In this section we describe solutions of the functional equation
\be \la{gopez} h=f\circ p=g\circ q,\ee where $f:\, C_1\rightarrow \C\P^1,$ $g:\, C_2\rightarrow \C\P^1$ are fixed 
holomorphic
functions on fixed Riemann surfaces $C_1,C_2$ and $h:\, C\rightarrow \C\P^1,$
$p:\, C\rightarrow C_1,$ $q:\, C\rightarrow C_2$ are 
unknown holomorphic functions on an unknown Riemann surface $C$. We always will assume that the considered Riemann surfaces 
are connected and compact.
\subsection{\la{prlm} Preliminaries}
Let $S\subset \C\P^1$ be a finite set 
and $z_0$ be a point from \linebreak $\C\P^1\setminus S$.
Recall that 
for any collection consisting of a Riemann surface $R$, holomorphic function 
$p:\, R\rightarrow \C\P^1$ non ramified outside 
of $S$, and a point $e\in p^{-1}\{z_0\}$ the homomorphism of the fundamental groups 
$$p_{\star}:\ \pi_1(R\setminus p^{-1}\{S\},e)\rightarrow 
\pi_1(\C\P^1\setminus S,z_0)$$ is a monomorphism such that its
image $\Gamma_{p,e}$ is a subgroup of finite index in the group $\pi_1(\C\P^1\setminus S,z_0)$,
and 
vice versa if $\Gamma$ is  
a subgroup of finite index in 
$\pi_1(\C\P^1\setminus S,z_0)$ then
there exist a Riemann surface $R$,
a function  $p:\, R\rightarrow \C\P^1$, and a point $e\in p^{-1}\{z_0\}$ 
such that $$p_{\star}(\pi_1(R\setminus p^{-1}\{S\},e))=\Gamma.$$ 
Furthermore, this correspondence descends to a one-to-one correspondence 
between conjugacy classes of subgroup of index $d$ in 
$\pi_1(\C\P^1\setminus S,z_0)$
and equivalence classes of holomorphic functions of degree $d$
non ramified outside of $S$, where 
functions $p:\, R\rightarrow \C\P^1$ and  $\tilde p:\, \tilde R\rightarrow \C\P^1$ are
considered as equivalent if there exists an isomorphism $w:\, R\rightarrow \tilde R $ 
such that $p=\tilde p\circ w.$

For collections $p_1:\, R_1\rightarrow \C\P^1,$ $e_1\in p_1^{-1}\{z_0\}$ and  
$p_2:\, R_2\rightarrow \C\P^1,$ $e_2\in p_2^{-1}\{z_0\}$
the groups $\Gamma_{p_1,e_1}$ and $\Gamma_{p_2,e_2}$ coincide if and only if there exists 
an isomorphism $w:\, R_1\rightarrow  R_2$ 
such that $p_1=p_2\circ w$ and $w(e_1)=e_2.$ 
More generally, the inclusion $$\Gamma_{p_1,e_1}\subseteq \Gamma_{p_2,e_2}$$
holds if and only if there exists a holomorphic function $w:\ R_1 \rightarrow R_2$ such that 
$p_1=p_2\circ w$ and 
$w(e_1)=e_2$ and
in the case if such a function exists it is defined in a unique 
way. Notice that this implies 
that if $p:\, R\rightarrow \C\P^1$, $e\in p^{-1}\{z_0\}$ is a pair
such that 
\be \la{cyki} \Gamma_{p_1,e_1}\subseteq \Gamma_{p,e} \subseteq \Gamma_{p_2,e_2}\ee
and $v:\, R_1\rightarrow R,$ $u:\, R\rightarrow R_2,$ are holomorphic function such that 
$p=p_2\circ u,$ $p_1=p\circ v$ and $v(e_1)=e,$ $u(e)=e_2$ 
then $w=u\circ v$. In particular, the function $w$ can be decomposed into a composition of 
holomorphic functions of degree greater than 1 if and only if there exists $\Gamma_{p,e}$ 
distinct from $\Gamma_{p_1,e_1}$ and $\Gamma_{p_2,e_2}$ such \eqref{cyki} holds.

In view of the fact that holomorphic 
functions can be identified with coverings of Riemann surfaces all the results above follow from the corresponding results about coverings (see e.g. 
\cite{mas}). Notice that the more customary language describing compositions of coverings uses monodromy groups of the functions involved rather than subgroups of $\pi_1(\C\P^1\setminus S,z_0)$. The interaction between these languages is explained below. 
In the paper we will use both these languages.

Fix a numeration $\{z_1, z_2, \dots , z_r\}$ of points of $S$ and for
each $i,$ $1\leq i \leq r,$ fix
a small loop $\beta_i$ around $z_i$ so that $\beta_1\beta_2 \dots \beta_r=1$ in $\pi_1(\C\P^1\setminus S,z_0)$. If $p:\ R\rightarrow \C\P^1$ 
is a holomorphic function non ramified outside of $S$ then for
each $i,$ $1\leq i \leq r,$ the loop
$\beta_i$ after the lifting by $p$ induces a permutation $\alpha_i(p)$ of points of $p^{-1}\{z_0\}$. The group $G_{p}$ generated by $\alpha_i(p),$ $1\leq i \leq r,$ is called the monodromy group of $p.$ Clearly, the group $G_{p}$ 
is transitive and the equality $\alpha_1(p)\alpha_2(p) \dots \alpha_r(p)=1$ 
holds in $G_{p}$. The representation of $\alpha_i(p),$ $1\leq i \leq r,$ 
by elements of the corresponding symmetric group
depends on the numeration of points of $p^{-1}\{z_0\}$
but the  
conjugacy class of the corresponding collection of permutations 
is well defined. Moreover, there is a one-to-one correspondence between
equivalence classes of holomorphic functions of degree $d$ non ramified outside of $S$ and 
conjugacy classes of ordered collections   
of permutations $\alpha_i,$ $1\leq i \leq r,$ from 
the symmetric group $S_d$ acting on the set $\{1,2,\dots, d\}$
such that $\alpha_1\alpha_2 
\dots \alpha_r=1$ and the permutation group generated by $\alpha_i,$ $1\leq i \leq r,$ is transitive (see e.g. \cite{mir}, Corollary 4.10). We will denote the conjugacy  
class of permutations which corresponds to a holomorphic function $p:\, R\rightarrow \C\P^1$
by $\hat\alpha(p).$ If $$ \phi_p:\, \pi_1(\C\P^1\setminus S,z_0)\rightarrow G_p\subset S_d$$ is a homomorphism which sends
$\beta_i$ to $\alpha_i$, $1\leq i \leq r,$ then the set of preimages of the stabilizers $G_{p,i},$ $1\leq i \leq d,$ coincides with the set of the groups $\Gamma_{p,e},$ $e\in p^{-1}\{z_0\}.$ 
On the other hand, for any group $\Gamma_{p,e}$, $e\in p^{-1}\{z_0\}$ the collection
of permutations $\alpha_i,$ $1\leq i \leq r,$ induced on the cosets of $\Gamma_{p,e}$ by $\beta_i,$ 
$1\leq i \leq r,$ is 
a representative of $\hat\alpha(p).$

If a holomorphic function 
$p:\, R\rightarrow \C\P^1$ of degree $d$ can be decomposed into a composition $p=f\circ q$ of holomorphic 
functions $q:\, R\rightarrow C$ and $f:\, C\rightarrow \C\P^1$ 
then the group $G_{p}$ has an imprimitivity system $\Omega_f$ consisting of 
$d_1=\deg f$ blocks 
such that the collection of permutations of 
blocks of $\Omega_f$ 
induced by $\alpha_i(p),$ $1\leq i \leq r,$
is a representative of 
$\hat\alpha(f)$, and vice versa if $G_{p}$ has an imprimitivity system $\Omega$ such that the collection of 
permutations of blocks of $\Omega$ 
induced by $\alpha_i(p),$ $1\leq i \leq r,$ is a representative of 
$\hat\alpha(f)$ for some holomorphic 
function $f:\, C\rightarrow \C\P^1$
then there exists a function $q:\, R\rightarrow C$ such that $p=f\circ q$. 
Notice that if the set $\{1,2,\dots, d\}$ is identified with the set $p^{-1}\{z_0\}$
then the set of blocks of the imprimitivity system $\Omega_f$ corresponding to the decomposition  
$p=f\circ q$  
has the form 
$\f B_i=q^{-1}\{t_i\},$ $1\leq i \leq d_1,$ where $\{t_1,t_2,\dots, t_{d_1}\} =f^{-1}\{z_0\}$.

If $p=\tilde f\circ \tilde q$, where $\tilde f:\, \tilde C\rightarrow \C\P^1$, 
$\tilde q:\, R\rightarrow \tilde C,$
is an other decomposition of $p$ then  
the imprimitivity systems $\Omega_f$, $\Omega_{\tilde f}$ coincide if and only there exists an automorphism $\mu:\, \tilde C \rightarrow C$
such that $$ f=\tilde f \circ \mu^{-1}, \ \ \ q=\mu \circ \tilde q.$$
In this case the decompositions $f\circ q$ and $\tilde f\circ \tilde q$
are called equivalent.
Therefore, equivalence classes of decompositions of $p$ are in a one-to-one correspondence
with imprimitivity systems of $G_p.$ More generally, if
$\f B$ is a block of $\Omega_{f}$ and $\f C$ is a block of $\Omega_{\tilde f}$ 
such that $\f B\cap \f C$ is non-empty, then $\f B$ and $\f C$    
have an intersection of cardinality $l$ if and only if there exist 
holomorphic functions $w:\, R \rightarrow R_1,$ $q_1:\, R_1\rightarrow C,$
$\tilde q_1:\, R_1\rightarrow \tilde C$, where $\deg w=l$, such that 
$$q=q_1\circ w, \ \ \ \ \ \ \tilde q=\tilde q_1\circ w.$$
In particular, if $p=f\circ q =f \circ q_1$ and the imprimitivity systems corresponding to 
the decompositions $p=f\circ q$ and $p=f\circ q_1$ coincide then $q_1 =\omega \circ q$ 
where $\omega$ is an automorphism of the surface $C$ 
such that $f \circ \omega=f$. Notice however that in general the equality $f\circ q=f\circ q_1$ 
does not imply that $q_1 =\omega \circ q$ for some $\omega$ as above. On the other hand,
since a holomorphic function $q:\, R\rightarrow C$ takes all the values on $C$ the equality 
$f\circ q=f_1\circ q$ always implies that $f=f_1.$

By the analogy with rational functions we will call a holomorphic function \linebreak
$p:\, R\rightarrow \C\P^1$ of degree greater than 1 indecomposable if the equality $p=f\circ q$ for some holomorphic 
functions $q:\, R\rightarrow C$ and $f:\, C\rightarrow \C\P^1$ implies that at least one of the functions 
$f,q$ is of degree 1. 
Clearly, if $p$ is non-ramified outside 
of $S$ and $z_0\in \C\P^1\setminus S$ then $p$ is 
indecomposable if and only if the subgroups $\Gamma_{p,e},$ $e\in p^{-1}\{z_0\}$ are 
maximal in $\pi_1(\C\P^1\setminus S,z_0).$


\subsection{\la{fsys} 
Description of solutions of equation \eqref{gopez}} 
Let $S=\{z_1, z_2, \dots , z_r\}$ be a union of branch points
of $f,g$ and $z_0$ be a fixed point from $\C\P^1\setminus S$. 

\bp \la{p1} Let $f:\, C_1\rightarrow \C\P^1,$ $g:\, C_2\rightarrow \C\P^1$ be 
holomorphic functions.
Then for any $a\in f^{-1}\{z_0\}$ and 
$b\in g^{-1}\{z_0\}$ there exist
holomorphic functions $u:\, C\rightarrow C_1,$ $v:\, C\rightarrow C_2$,
$h:\, C\rightarrow \C\P^1,$ and a point $c\in h^{-1}\{z_0\}$ 
such that 
\be \la{f1} h=f\circ u=g\circ v, \ \ \ u(c)=a, \ \ \ v(c)=b.\ee
Furthermore, the function $h$ has the following property: if 
\be \la{f2} \tilde h=f\circ \tilde u=g\circ \tilde v, \ \ \  \tilde u(\tilde c)=a, 
\ \ \ \tilde v(\tilde c)=b \ee
for some
holomorphic functions $\tilde h:\, \tilde C\rightarrow \C\P^1,$ $\tilde u:\, \tilde C\rightarrow C_1,$ $\tilde v:\, \tilde C\rightarrow C_2$,
and a point $\tilde c\in \tilde h^{-1}\{z_0\},$ then
there exists a holomorphic function $w:\, \tilde C\rightarrow C$ such that
\be \la{e3} \tilde h= h\circ w, \ \ \ \tilde u= u\circ  w,\ \ \ \tilde v= v\circ w,
\ \ \ w(\tilde c)=c.\ee
\ep
\pr Since the subgroups $\Gamma_{f,a}$ and $\Gamma_{g,b}$ are of finite index in 
$\pi_1(\C\P^1\setminus S,z_0)$ their intersection is also of finite index. Therefore, there exists 
a pair $h:\, C\rightarrow \C\P^1,$ $c\in h^{-1}\{z_0\}$ such that 
$\Gamma_{h,c}=\Gamma_{f,a}\cap\Gamma_{g,b}$ and for such a pair equalities \eqref{f1} hold. 
Furthermore, equalities \eqref{f2} imply that 
$\Gamma_{\tilde h,\tilde c}\subseteq \Gamma_{f,a}\cap\Gamma_{g,b}$
Therefore, 
$\Gamma_{\tilde h,\tilde c}\subseteq \Gamma_{h,c}$ and 
hence 
$\tilde h=h\circ w$
for some $w:\, \tilde C\rightarrow C$  
such that $w(\tilde c)=c.$ 
It follows now from $$f\circ \tilde u=f\circ u \circ w, \ \ \ g\circ \tilde v=g\circ v \circ w$$ and 
$$(u\circ w)(\tilde c)=\tilde u(\tilde c), \ \ \ \ (v\circ w)(\tilde c)=\tilde v(\tilde c)$$ 
that
$$\tilde u=u\circ w, \ \  \  \ \tilde v=v\circ w.\ \ \ \ \ \Box$$    
\vskip 0.2cm

\vskip 0.2cm
For holomorphic functions $f:\, C_1\rightarrow \C\P^1,$ $\deg f=n,$ and $g:\, C_2\rightarrow \C\P^1,$ $\deg g=m,$ 
fix some representatives
$\alpha_i(f),$ $\alpha_i(g),$ $1\leq i \leq r,$ of the classes
$\hat\alpha(f),$ $\hat\alpha(g)$ and 
define the permutations $\delta_1, \delta_2, \dots, \delta_{r}\in S_{nm}$ on
the set of $mn$ elements $c_{j_1,j_2},$ $1\leq j_1 \leq n,$ $1\leq j_2\leq m,$ 
as follows: $c_{j_1,j_2}^{\delta_i}=c_{j_1^{\prime},j_2^{\prime}},$ where $$j_1^{\prime}
=j_1^{\alpha_i(f)},\ \ \ \ j_2^{\prime}=j_2^{\alpha_i(g)}, \ \ \ \ 1\leq i \leq r.$$ 
It is convenient to consider $c_{j_1,j_2},$ $1\leq j_1 \leq n,$
$1\leq j_2\leq m,$ as elements of a $n\times m$  matrix $M$. Then the action of the 
permutation $\delta_i,$ $1\leq i \leq r,$
reduces to the permutation of rows of $M$ in accordance with the permutation $\alpha_i(f)$ and  
the permutation of columns of $M$ in accordance with the permutation $\alpha_i(g)$.

In general 
the permutation group $\Gamma(f,g)$ generated by $\delta_i,$ $1\leq i \leq r,$
is not transitive on 
the set $c_{j_1,j_2},$ $1\leq j_1 \leq n,$
$1\leq j_2\leq m$. Denote by $o(f,g)$ the number of transitivity sets of 
$\Gamma(f,g)$ and let $\delta_{i}(j),$ $1\leq j \leq o(f,g),$
$1\leq i \leq r,$ be the 
permutation induced by the permutation $\delta_i,$ $1\leq i \leq r,$
on the transitivity set $U_j,$ $1\leq j \leq o(f,g).$ 
By construction, for any $j,$ 
$1\leq j \leq o(f,g),$ the permutation group $G_j$ generated by $\delta_i(j),$ 
$1\leq i \leq r,$ is transitive and the
equality 
$$\delta_{1}(j)\delta_{2}(j)\dots \delta_{r}(j)=1$$
holds. Therefore, there exist holomorphic functions 
$h_j:\, R_j\rightarrow \C\P^1,$ $1\leq j \leq o(f,g),$ 
such that the collection $\delta_{i}(j),$ 
$1\leq i \leq r,$ is a representative of $\hat\alpha(h_j).$
Moreover, it follows from the construction that for each $j,$ $1\leq j \leq o(f,g),$
the intersections of the transitivity set $U_j$ with rows of $M$
form an imprimitivity system $\Omega_f(j)$ for $G_j$ such that the permutations
of blocks of $\Omega_f(j)$ induced by 
$\delta_i(j),$ $1\leq i \leq r,$ coincide with $\alpha_i(f).$ 
Similarly, the intersections of $U_j$ with columns of $M$
form an imprimitivity system $\Omega_g(j)$ such that the permutations
of blocks of $\Omega_g(j)$ induced by 
$\delta_i(j),$ $1\leq i \leq r,$ coincide with $\alpha_i(g).$ This implies that there exist holomorphic
functions $u_j:\,  R_j\rightarrow C_1$ and $v_j:\,  R_j\rightarrow C_2$
such that 
\be \la{e4} h_j=f\circ u_j=g\circ v_j.\ee

\bt \la{p2} Let $f:\, C_1\rightarrow \C\P^1,$ $g:\, C_2\rightarrow \C\P^1$ be 
holomorphic functions. 
Suppose that $h:\, R\rightarrow \C\P^1,$ $p:\, R\rightarrow C_1,$ $q:\, R\rightarrow C_2$
are holomorphic 
function such that \be  \la{ur} h=f\circ p=g\circ q.\ee Then
there exist $j,$ $1\leq j \leq o(f,g),$ and holomorphic functions $w:\, R\rightarrow R_j,$ 
$\tilde p:\, R_j\rightarrow C_1,$ $\tilde q:\, R_j\rightarrow C_2$ 
such that 
\be \la{e5} h= h_j\circ w, \ \ \ p= \tilde p\circ w,\ \ \ q= \tilde q\circ w
\ee
and
$$f\circ \tilde p \sim f \circ u_j,\ \ \ \ \ \ g\circ \tilde 
q\sim g\circ v_j.$$
\et 
 
\pr 
It follows from Proposition \ref{p1} that in order to prove the theorem it 
is enough to show that for any choice of points $a\in f^{-1}\{z_0\}$ and 
$b\in g^{-1}\{z_0\}$ the class of permutations $\hat\alpha(h)$ corresponding to the function $h$  
from Proposition \ref{p1}
coincides with $\hat\alpha(h_j)$ for some $j,$ $1\leq j \leq o(f,g).$
On the other hand, the last statement is equivalent to the statement that 
for any choice $a\in f^{-1}\{z_0\}$ and  
$b\in g^{-1}\{z_0\}$  
there exists $j,$ $1\leq j \leq o(f,g),$ and an element $c$
of the transitivity set $U_j$ such that the group $\Gamma_{f,a}\cap\, \Gamma_{g,b}$ 
is the preimage of the stabilizer $G_{j,c}$ of $c$ in the group $G_j$ under the homomorphism 
$$\phi_{h_j}:\, \pi_1(\C\P^1\setminus S,z_0)\rightarrow G_j$$ (see Subsection \ref{prlm}).

For fixed $a\in f^{-1}\{z_0\},$
$b\in g^{-1}\{z_0\}$ 
let $l$ be the index which corresponds to the point $a$  
under the identification of the set $f^{-1}\{z_0\}$ with the set
$\{1,2, \dots,n\}$, $k$ be the index which corresponds to the point $b$ 
under the identification of the set $g^{-1}\{z_0\}$
with the set $\{1,2, \dots,m\}$, and $U_{j}$ be the transitivity set of 
$\Gamma(f,g)$ containing the element $c_{l,k}.$ 
We have: \be \la{trr} \Gamma_{f,a}=\phi_f^{-1}\{G_{f,\,l}\}, \ \ \ \ \Gamma_{g,b}=\phi_g^{-1}\{G_{g,k}\}.\ee
Furthermore, if $\psi_1:\, G_f\rightarrow G_j$ (resp. $\psi_2:\, G_g\rightarrow G_j$) is a homomorphism which sends
$\alpha_i(f)$ (resp. $\alpha_i(g)$) to  
$\alpha_i(h_j),$ $1\leq i \leq r,$ then 
\be \la{trr1} G_{f,\,l}=\psi_1^{-1}\{A_l\}, \ \ \ \ \  G_{g,k}=\psi_2^{-1}\{B_k\},\ee
where $A_l$ (resp. $B_k$) is the subgroup of $G_j$ which transforms the set of elements $c_{j_1,j_2}\in U_j$
for which $j_1=a$ (resp. $j_2=b$) to itself. 

Since $$ \psi_1\circ \phi_f=\psi_2\circ \phi_g=\phi_{h_j}$$ it follows from \eqref{trr}, \eqref{trr1} that
$$\Gamma_{f,a}\cap \Gamma_{g,b}= (\psi_1\circ \phi_f)^{-1}\{A_l\}\cap (\psi_2\circ \phi_g)^{-1}\{B_k\}=$$ 

$$=
\phi_{h_j}^{-1}\{A_l\}\cap \phi_{h_j}^{-1}\{B_k\}=
\phi_{h_j}^{-1}\{A_l \cap B_k\}=\phi_{h_j}^{-1}\{G_{j,c_{k,l}}\}. \ \ \ \ \ \Box$$

\vskip 0.2cm

For $i,$ $1\leq i \leq r,$ denote by
$$
\l_i=(f_{i,1},f_{i,2}, ... , f_{i,u_i})
$$
the collection of lengths of disjoint cycles in the permutation $\alpha_i(f)$,
by 
$$
\mu_i=(g_{i,1},g_{i,2}, ... , g_{i,v_i}) 
$$
the collection of lengths of disjoint cycles in the permutation $\alpha_i(g),$ and 
by $g(R_j)$, $1\leq j \leq o(f,g),$ the genus of the surface $R_j$.
The proposition 
below generalizes the corresponding result of Fried (see \cite{f3}, Proposition 2)
concerning the case when $f,g$ are rational functions.

\bp \la{p3} In the above notation the formula 
\be \la{rh0} 
\sum_{j=1}^{o(f,g)}(2-2g(R_j)) =
\sum_{i=1}^{r}
\sum_{j_1=1}^{u_{i}} \sum_{j_2=1}^{v_{i}} \GCD(f_{i,j_1}g_{i,j_2})-(r-2)nm
\ee
holds.
\ep
\pr
Denote by $e_i(j),$ $1\leq i \leq r,$ $1\leq j \leq o(f,g),$
the number of disjoint cycles in the permutation 
$\delta_i(j)$. Since for any $j,$ $1\leq j \leq o(f,g),$ the Riemann-Hurwitz formula implies that 
$$2-2g(R_j)=
\sum_{i=1}^{r}e_{i}(j)-(r-2)
\vert U_j\vert$$ we have: 
$$\sum_{j=1}^{o(f,g)}(2-2g(R_j))=\sum_{j=1}^{o(f,g)}\sum_{i=1}^
{r}e_{i}(j)-(r-2)mn.$$

On the other hand, it follows from the construction 
that for given $i,$ $1\leq i \leq r,$ 
$$\sum_{j=1}^{o(f,g)}e_{i}(j)=\sum_{j_1=1}^{u_{i}}\sum_{j_2=1}^{v_{i}} \GCD(f_{i,j_1}g_{i,j_2})$$ 
and hence 
$$\sum_{j=1}^{o(f,g)}\sum_{i=1}^
{r}e_{i}(j)=\sum_{i=1}^{r}
\sum_{j_1=1}^{u_{i}} \sum_{j_2=1}^{v_{i}} \GCD(f_{i,j_1}g_{i,j_2}). \ \ \ 
\Box$$
\vskip 0.2cm

The proposition below shows that if $f,g$ are rational functions then the Riemann surfaces $R_j,$ $1\leq j \leq o(f,g),$ 
may be identified with irreducible components of the 
affine algebraic curve 
$$ h_{f,g}(x,y):\ P_1(x)Q_2(y)-P_2(x)Q_1(y)=0,$$ 
where 
$P_1$, $P_2$ and $Q_1,Q_2$ are pairs polynomials without common roots such that 
$$ f=P_1/P_2, \ \ g=Q_1/Q_2.$$

\bp \la{compon} For rational functions $f,g$ the corresponding Riemann surfaces $R_j,$  $1\leq j \leq o(f,g),$
are in a one-to-one correspondence with irreducible components of the curve $h_{f,g}(x,y).$
Furthermore, each $R_j$ is a desingularization of the corresponding component.
In particular, the curve $h_{f,g}(x,y)$ is irreducible if and only if the group $\Gamma(f,g)$ is transitive.
\ep

\pr For $j,$ $1\leq j \leq o(f,g),$ 
denote by $S_j$ the union of poles of $u_j$ and $v_j$ and 
define the mapping $t_j:\, R_j\setminus S_j \rightarrow \C^2$  
by the formula $$z\rightarrow (u_j,v_j).$$
It follows from formula \eqref{e4} that for each $j,$ $1\leq j \leq o(f,g),$ the mapping 
$t_j$ 
maps $R_j$ to an irreducible component of the curve $h_{f,g}(x,y)$. Furthermore, for any 
point $(a,b)$ on $h_{f,g}(x,y)$, such that $z_0=f(a)=g(b)$ is not contained in $S$, 
there exist uniquely defined $j,$ $1\leq j \leq o(f,g),$ and $c\in h_j^{-1}\{z_0\}$
satisfying $$ u_j(c)=a, \ \ \ v_j(c)=b.$$
This implies that the Riemann surfaces $R_j,$ $1\leq j \leq o(f,g),$
are in a one-to-one correspondence with irreducible components of $h_{f,g}(x,y)$ 
and that each mapping $t_j,$ $1\leq j \leq o(f,g),$
is generically injective. Since an injective mapping of Riemann surfaces is an isomorphism onto an open subset 
we conclude that 
each $R_j$ is a desingularization of the corresponding component of $h_{f,g}(x,y)$. \qed

\section{\la{s2} Irreducible and reducible pairs}
Let $f:\, C_1\rightarrow \C\P^1,$ 
$g:\, C_2\rightarrow \C\P^1$ be a pair of holomorphic functions
non-ramified outside of $S$ and $z_0\in \C\P^1\setminus S$. 
By the analogy with the rational case we will call the pair $f,g$ 
{\it irreducible} if $o(f,g)=1$.
Otherwise we will call such the pair $f,g$
{\it reducible}.  
In this section we study properties of irreducible and reducible pairs.

\bp \la{p4} A pair of holomorphic functions $f:\, C_1\rightarrow \C\P^1,$ 
$g:\, C_2\rightarrow \C\P^1$ is irreducible whenever their degrees are coprime.
\ep
\pr Let $n=\deg f, $ $m=\deg g.$
Since the index of $\Gamma_{f,a}\cap \Gamma_{g,b}$ coincides with the 
cardinality of the corresponding imprimitivity set $U_j,$ the pair $f,g$ is irreducible if and only if 
for any $a\in f^{-1}\{z_0\},$ $b\in g^{-1}\{z_0\}$ the equality 
\be \la{z} [\pi_1(\C\P^1\setminus S,z_0):\Gamma_{f,a}\cap \Gamma_{g,b}]=nm\ee
holds. Since the index of $\Gamma_{f,a}\cap \Gamma_{g,b}$ in 
$\pi_1(\C\P^1\setminus S,z_0)$ is a multiple of the indices of
$\Gamma_{f,a}$ and $\Gamma_{g,b}$ in $\pi_1(\C\P^1\setminus S,z_0)$,
this index is necessary equal to $mn$ whenever $n$ and $m$ are coprime.  \qed 

\bt \la{p5}  A pair of holomorphic functions $f:\, C_1\rightarrow \C\P^1,$ 
$g:\, C_2\rightarrow \C\P^1$ is irreducible if and only if  
for any $a\in f^{-1}\{z_0\},$ 
$b\in g^{-1}\{z_0\}$ 
the equality 
\be \la{zz} \Gamma_{f,a}\Gamma_{g,b}=\Gamma_{g,b}\Gamma_{f,a}=\pi_1(\C\P^1\setminus S,z_0)
\ee
holds. 
\et

\pr
Since $$\left[\pi_1(\C\P^1\setminus S,z_0)
:\Gamma_{f,a}\cap \Gamma_{g,b}\right]
=
\left[\pi_1(\C\P^1\setminus S,z_0):\Gamma_{g,b}\right]\left[\Gamma_{g,b}^{\,}
:\Gamma_{f,a}\cap \Gamma_{g,b}\right],$$ the 
equality \eqref{z} is equivalent to the equality 
\be \la{z3} \left[\Gamma_{g,b}:\Gamma_{f,a}\cap \Gamma_{g,b}\right]=n.\ee
Recall that for any subgroups $A,B$ of finite index in a group $G$ the inequality 
\be \la{ggh} \left[<A,B>:A\right]\geq \left[B:A\cap B\right] \ee  
holds and the equality attains if and only if the groups $A$ and $B$ are permutable
(see e.g. \cite{kur}, p. 79). Therefore, 
$$n=\left[\pi_1(\C\P^1\setminus S,z_0):\Gamma_{f,a}\right]\geq 
\left[<\Gamma_{f,a},\Gamma_{g,b}
>:\Gamma_{f,a}\right]\geq \left[\Gamma_{g,b}:\Gamma_{f,a}\cap \Gamma_{g,b}\right]$$
and hence equality \eqref{z3} holds if and only if 
$\Gamma_{f,a}$ and $\Gamma_{g,b}$ are permutable and generate $\pi_1(\C\P^1\setminus S,z_0)$. 
\qed

\bc Let $f:\, C_1\rightarrow \C\P^1,$ 
$g:\, C_2\rightarrow \C\P^1$ be an irreducible
pair of holomorphic functions. Then any pair of holomorphic functions $\tilde f:\, \tilde C_1\rightarrow \C\P^1,$ 
$\tilde g:\, \tilde C_2\rightarrow \C\P^1$ such that 
$$f=\tilde f \circ p, 
\ \ \ \ g=\tilde g \circ q$$ for some 
holomorphic functions $p:\, C_1\rightarrow \tilde C_1,$ 
$q:\, C_2\rightarrow \tilde C_2$ is also irreducible.
\ec

\pr Since for any $\tilde a\in \tilde f^{-1}\{z_0\},$ 
$\tilde b\in \tilde g^{-1}\{z_0\}$ and $a\in p^{-1}\{\tilde a\},$ 
$b\in q^{-1}\{\tilde b\}$
the inclusions 
$$\Gamma_{f,a} \subseteq \Gamma_{\tilde f, \tilde a}, \ \ \ \ \Gamma_{g,b} \subseteq \Gamma_{\tilde g, \tilde b}$$ hold
it follows from \eqref{zz} that $$\Gamma_{\tilde f, \tilde a}\Gamma_{\tilde g, \tilde b}=\Gamma_{\tilde g, \tilde b}
\Gamma_{\tilde f, \tilde a}=
\pi_1(\C\P^1\setminus S,z_0).\ \ \ \ \Box$$ 

\vskip 0.2cm


Set $$ \Gamma_{N_g}=\bigcap_{b\in g^{-1}\{z_0\}} \Gamma_{g,b}$$
and denote by $\hat{N}_g$
the corresponding equivalence class of holomorphic functions.
Since the subgroup $ \Gamma_{N_g}$ is normal in $\pi_1(\C\P^1\setminus S,z_0)$, for any 
$a_1, a_2\in f^{-1}\{z_0\}$ 
the subgroups $\Gamma_{f,a_1}\Gamma_{N_g}$ and $\Gamma_{f,a_2}\Gamma_{N_g}$
are conjugated. We will denote the equivalence class of holomorphic functions 
corresponding to this conjugacy class by $f\hat {N}_g$. 

\bp \la{kij} For any pair of holomorphic functions $f:\, C_1\rightarrow \C\P^1,$ \linebreak $g:\, C_2\rightarrow \C\P^1$ 
and a representative $fN_g:\, C\rightarrow \C\P^1$ of $f\hat {N}_g$
the equality 
$$o(f,g)=o(fN_g,g)$$ holds. 
\ep
  
\pr 
For any $a\in f^{-1}\{z_0\},$ $b\in g^{-1}\{z_0\}$ the action of the permutation group $\Gamma(f,g)$ can be identified with the action of 
$\pi_1(\C\P^1\setminus S,z_0)$ on pairs of cosets 
$\alpha_{j_1}\Gamma_{f,a},$ $\beta_{j_2}\Gamma_{g,b},$ $1\leq j_1 \leq n,$ $1\leq j_2 \leq m.$ 
Furthermore, two pairs $\alpha_{j_1}\Gamma_{f,a},$ $\beta_{j_2}\Gamma_{g,b}$ and $\alpha_{i_1}\Gamma_{f,a},$ $\beta_{i_2}\Gamma_{g,b}$ are in the same orbit
if and only if the set \be \la{cfr} \alpha_{i_1}\Gamma_{f,a}\alpha_{j_1}^{-1}\cap \beta_{i_2}\Gamma_{g,b}\beta_{j_2}^{-1}\ee is non-empty.

Associate now to an orbit $\Gamma(f,g)$ containing the pair $\alpha_{j_1}\Gamma_{f,a},$ $\beta_{j_2}\Gamma_{g,b},$
$1\leq j_1 \leq n,$ $1\leq j_2 \leq m,$
an orbit of $\Gamma(fN_g,g)$ containing the pair 
$\alpha_{j_1}\Gamma_{f,a}\Gamma_{N_g},$ $\beta_{j_2}\Gamma_{g,b}.$
If set \eqref{cfr} is non-empty then the set \be\la{cfr1}   \alpha_{i_1}\Gamma_{f,a}\Gamma_{N_g}\alpha_{j_1}^{-1}\cap \beta_{i_2}\Gamma_{g,b}\beta_{j_2}^{-1}\ee 
is also non-empty and therefore we obtain a well-defined
map $\phi$ from the set of orbits of $\Gamma(f,g)$ to the set of orbits of $\Gamma(fN_g,g)$. Besides, the map $\phi$ is clearly surjective. 

In order to prove the injectivity of $\phi$ we must show that if set \eqref{cfr1} is non-empty then set \eqref{cfr} is also non-empty. So suppose that \eqref{cfr1} is non-empty and let
$x$ be its element. 
In view of the normality of $\Gamma_{N_g}$ the equality 
$$  \alpha_{i_1}\Gamma_{f,a}\Gamma_{N_g}\alpha_{j_1}^{-1}=  \alpha_{i_1}\Gamma_{f,a}\alpha_{j_1}^{-1}\Gamma_{N_g}$$  
holds and therefore
there exist $\alpha\in \Gamma_{f,a}$, $\beta \in\Gamma_{N_g}$, and $\gamma\in \Gamma_{g,b}$ such that 
$$x=\alpha_{i_1}\alpha \alpha_{j_1}^{-1}\beta  =\beta_{i_2}\gamma \beta_{j_2}^{-1}.$$ Furthermore, it follows from the definition of 
$\Gamma_{N_g}$ that there exists $\gamma_1\in \Gamma_{g,b}$ such that $\beta = \beta_{j_2}\gamma_1 \beta_{j_2}^{-1}$.
Set $y=x\beta^{-1}$. Then we have: $$y=\alpha_{i_1}\alpha \alpha_{j_1}^{-1}=\beta_{i_2}\gamma \beta_{j_2}^{-1}\beta^{-1}=\beta_{i_2}\gamma\gamma_1^{-1} \beta_{j_2}^{-1}.$$ This implies that $y$ 
is contained in set
\eqref{cfr} and hence \eqref{cfr} is non-empty. \qed

The following results generalizes the corresponding result of Fried about rational
functions (see \cite{f2}, Proposition 2).

\bt \la{p6} For any reducible pair of holomorphic functions $f:\, C_1\rightarrow \C\P^1,$ $g:\, C_2\rightarrow \C\P^1$ there exist
holomorphic functions  $f_1:\, \tilde C_{1}\rightarrow \C\P^1,$ $g_1:\,\tilde C_{2}\rightarrow \C\P^1,$
and $p:\, C_{1}\rightarrow \tilde C_{1},$ $q:\, C_{2}\rightarrow \tilde C_{2}$
such that
\be \la{e6} f= f_1\circ p, \ \ \ g=g_1\circ q,\ \ \ o(f,g)=o (f_1, g_1),\qquad \mbox{and} \qquad \hat N_{f_1}=\hat N_{g_1}.\ee
\et 
 
\pr For a holomorphic function $p:\, R\rightarrow \C\P^1$ denote by  
$d(p)$ a maximal number such that there
exist holomorphic functions of degree greater than 1
$$p_1:\, R\rightarrow R_1, \ \ p_i:\, R_{i-1}\rightarrow R_i, \ \ 2\leq i \leq d(p)-1, \ \ p_{d(p)}:\, R_{d(p)-1}\rightarrow \C\P^1$$ 
such that $$p=p_{d(p)}\circ p_{d(p)-1}\circ \dots \circ p_{1}.$$
We use the induction on the number $d=d(f)+d(g)$.

If $d=2$ that is if both functions 
$f,g$ are indecomposable then the equality $d(f)=1$, taking into account the normality of $N_g$, implies 
that either \be \la{ry} \Gamma_{f,a}N_{g}=\Gamma_{f,a}\ee for all $a\in f^{-1}\{z_0\}$
or \be \la{ryry} \Gamma_{f,a}N_{g}=\pi_1(\C\P^1\setminus S,z_0)\ee
for all $a\in f^{-1}\{z_0\}.$
The last possibility however would imply that for any \linebreak $b\in g^{-1}\{z_0\}$
$$\Gamma_{f,a}\Gamma_{g,b}=\Gamma_{g,b}\Gamma_{f,a}=\pi_1(\C\P^1\setminus S,z_0)$$ in contradiction with Theorem \ref{p5}. Therefore,
equalities \eqref{ry} hold and hence
$$N_{g}\subseteq\bigcap_{a\in f^{-1}\{z_0\}} \Gamma_{f,a}= N_{f}.$$ The same arguments show that $N_{f}\subseteq N_{g}.$
Therefore, $N_{g}= N_{f}$ and we can set $f_1=f,$ $g_1=g.$ 

Suppose now that $d>2.$
If $N_{f}=N_{g}$ then as above we can set  
$f_1=f,$ $g_1=g$ so assume that $N_{f}\neq N_{g}$. Then, again taking into account the normality of $N_g$, either  
\be \la{xzc} \Gamma_{f,a}\subsetneq \Gamma_{f,a}N_{g},\ee for all $a\in f^{-1}\{z_0\}$ or 
$$\Gamma_{g,b}\subsetneq \ \Gamma_{g,b}N_{f}$$ for all $b\in g^{-1}\{z_0\}$.
Suppose say that \eqref{xzc} holds. Since equality \eqref{ryry} is impossible this implies that for any $a\in f^{-1}\{z_0\}$ there exist $h:\, C\rightarrow \C\P^1$ and $c\in h^{-1}\{z_0\}$ such that
$\Gamma_{f,a}N_{g}=\Gamma_{h,c}.$ 

It follows from \eqref{xzc} that 
$f=h\circ p$ for some $p:\, C_1\rightarrow C$ with $1<\deg h<\deg f$ and hence
$d(h)<d(f).$   
Since by Proposition \ref{kij} the equality $o(f,g)=o(h,g)$ holds
the theorem follows now from 
the induction assumption. \qed

\bt \la{rit2} Let $f:\, C_1\rightarrow \C\P^1,$ $g:\, C_2\rightarrow \C\P^1$ be an irreducible pair of holomorphic 
functions and $p:\, C\rightarrow C_1,$ $q:\, C\rightarrow C_2$ be holomorphic functions such that
$f\circ p=g\circ q.$ Suppose that $q$ is indecomposable. Then $f$ is also indecomposable.
\et 

\pr 
Set $h=f\circ p=g\circ q$ and fix a point $c\in h^{-1}\{z_0\}.$ 
Since
\be \la{incl} \Gamma_{h,c}\subseteq \Gamma_{f,a}, \ \ \ \
\Gamma_{h,c}\subseteq \Gamma_{g,b},\ee where $a=p(c), b=q(c),$ we have: 
\be \la{egik} \Gamma_{h,c}\subseteq \Gamma_{f,a}\cap \Gamma_{g,b}\subseteq  \Gamma_{g,b}.\ee
Furthermore, by Theorem \ref{p5} 
\be \la{32} \Gamma_{f,a} \Gamma_{g,b}=\pi_1(\C\P^1\setminus S,z_0).\ee 
Since \eqref{32} implies that $\Gamma_{f,a}\cap \Gamma_{g,b}\neq \Gamma_{g,b}$
it follows from \eqref{egik} taking into account 
the indecomposability of  
$q$ that 
\be \la{plm} \Gamma_{h,c}=\Gamma_{f,a}\cap \Gamma_{g,b}.\ee

In order to prove the theorem we must show that if $\Gamma\subseteq\pi_1(\C\P^1\setminus S,z_0)$ is a subgroup 
such that 
\be \la{fr} \Gamma_{f,a}\subsetneq\Gamma\ee 
then $\Gamma=\pi_1(\C\P^1\setminus S,z_0).$ Clearly, \eqref{32} implies that \be \la{322} \Gamma\, \Gamma_{g,b}=\pi_1(\C\P^1\setminus S,z_0).\ee 
Consider the intersection $$\Gamma_1=\Gamma\cap  \Gamma_{g,b}.$$ 
It follows from \eqref{ggh} and \eqref{32}, \eqref{322} that 
$$[\pi_1(\C\P^1\setminus S,z_0): \Gamma_{f,a}]= [\Gamma_{g,b}: \Gamma_{h,c}], \ \ \
[\pi_1(\C\P^1\setminus S,z_0): \Gamma]= [\Gamma_{g,b}: \Gamma_1].$$
Therefore, \eqref{fr} implies that $$[\Gamma_{g,b}: \Gamma_1]<[\Gamma_{g,b}: \Gamma_{h,c}]$$ and hence 
$\Gamma_{h,c}\subsetneq \Gamma_1.$ Since $\Gamma_1\subseteq \Gamma_{g,b}$ it follows now from the indecomposability of $q$ that $\Gamma_1= \Gamma_{g,b}.$ Therefore, $\Gamma_{g,b}\subseteq \Gamma.$ Since also $\Gamma_{f,a}\subseteq \Gamma$ it follows now from \eqref{32} that $\Gamma =\pi_1(\C\P^1\setminus S,z_0).$ \qed

\section{\la{s3} Double decompositions involving generalized polynomials}
Say that a holomorphic function $h:\, C\rightarrow \C\P^1$ is {\it a generalized polynomial} if $h^{-1}\{\infty\}$ 
consists of a unique point. In this section we mention some specific properties of double decompositions $f\circ p=g\circ q$ in the case when $f,g$ are generalized polynomials. 

We start from mentioning two corollaries of Theorem \ref{p6} for such double decompositions.

\bc \la{fr2} If in Theorem \ref{p6} the functions $f,g$ are generalized polynomials then $\deg f_1=
\deg g_1$.
\ec

\pr The equality $f=f_1\circ p$ for a generalized polynomial $f$ implies that $f_1$ is also a 
generalized polynomial.
Furthermore, since $\Gamma_{N_{f_1}}=\bigcap_{a\in f_1^{-1}\{z_0\}} \Gamma_{f_1,a}$ the monodromy group
of $\Gamma_{N_{f_1}}$ may be obtained by the repeated use of the construction given in Subsection \ref{fsys}.  
On the other hand, it is easy to see that if $f_1$ is a generalized polynomial 
then on each stage of this process the permutation 
corresponding to the loop around infinity consists of cycles
of length equal to the degree of $f_1$ only.  
Therefore, the same is true for $\Gamma_{N_{f_1}}$ and hence 
the equality $\hat N_{f_1}=\hat N_{g_1}$ implies that 
$\deg f_1=\deg g_1$. \qed
 
The following important specification of Theorem \ref{p6} goes back
to Fried (see \cite{f2}, Proposition 2). 

\bc \la{bbv} Let $A,B$ be polynomials such that curve \eqref{cur} is reducible. Then there exist polynomials $A_1,B_1, C,D$ such that 
\be \la{e66} A= A_1\circ C, \ \ \ B=B_1\circ D,\ \ \  \hat N_{A_1}=\hat N_{B_1},\ee
and each irreducible component $F(x,y)$ of curve \eqref{cur} has the form   
$F_1(C(x),D(y))$, where $F_1(x,y)$ is an irreducible component of the curve \be \la{kjk} A_1(x)-B_1(y)=0.\ee
\ec
\pr Indeed, it follows from Theorem \ref{p6} and Proposition \ref{compon} that there exist polynomials $A_1,B_1, C,D$ such that equalities
\eqref{e66} hold and curves \eqref{cur} and \eqref{kjk} have the same number of irreducible components. 
Since for each irreducible component $F_1(x,y)$ of curve \eqref{kjk} the polynomial $F_1(C(x),D(y))$ is a component of curve \eqref{cur} this implies that any irreducible component $F(x,y)$ of curve \eqref{cur} has the form $F_1(C(x),D(y))$ for some irreducible component $F_1(x,y)$ of curve \eqref{kjk}. \qed

For a holomorphic function $h:\, C\rightarrow \C\P^1$ and $z\in C$ denote by $\mm_z\, h$ the multiplicity of
$h$ at $z.$ 

\bt\la{polusa}
Let $f:\, C_1\rightarrow \C\P^1,$ $g:\, C_2\rightarrow \C\P^1$ be 
generalized polynomials, $\deg f=n,$ $\deg g=m,$ $l=\LCM(n,m),$
and $h:\, R\rightarrow \C\P^1,$ $p:\, R\rightarrow C_1,$ $q:\, R\rightarrow C_2$ be holomorphic 
functions such that 
\be \la{pl} h=f\circ p=g\circ q.\ee
Then there exist 
holomorphic functions $w:\, R\rightarrow C,$ 
$\tilde p:\, C\rightarrow C_1,$ $\tilde q:\, C\rightarrow C_2$ 
such that 
\be \la{ee} p= \tilde p\circ w,\ \ \ q= \tilde q\circ w,
\ee
and for any $z\in h^{-1}\{\infty\}$ 
$$\mm_z\, \tilde p=l/n, \ \ \ \mm_z\, \tilde q= l/m.$$\et

\pr In view of Theorem \ref{p2} it is enough to prove that if 
$u_j,v_j,$ $1\leq j \leq o(f,g),$
are functions defined in 
Subsection \ref{fsys} then for any $z\in h^{-1}\{\infty\}$ and $j,$ $1\leq j \leq o(f,g)$, the equalities
\be \la{iio} \mm_z\, u_j=l/n, \ \ \ \mm_z\, v_j=l/m\ee    
hold.

Since $f,g$ are generalized polynomials it follows from the construction given in
Subsection \ref{fsys} that for any function $h_j=f\circ u_j=g\circ v_j,$ $1\leq j \leq o(f,g),$ 
the permutation of its monodromy group corresponding to the loop around infinity consists of cycles
of length equal to $l$ only. On the other hand, the length of such a cycle coincides with the multiplicity
of the corresponding point from $h_j^{-1}\{\infty\}.$ Now equalities \eqref{iio} follow from the fact that 
for any $z\in R_j,$ $1\leq j \leq o(f,g),$
$$\mm_z h_j=\mm_{u_j(z)}f\,\mm_z u_j= \mm_{v_j(z)}g\, \mm_z v_j. \ \ \ \Box$$

\bc \la{copo} Let $A,$ $B$ be polynomials of the same degree $n$ and $C,$ $D$ be rational functions such that 
$$A\circ C=B\circ D.$$ Then there exist a
rational function $W$, mutually distinct points of the complex sphere $\gamma_i,$ $1\leq i \leq r,$ and complex numbers 
$\alpha_i, \beta_i$ $0\leq i \leq r,$
such that 
$$C=\left(\alpha_0+\frac{\alpha_1}{z-\gamma_1}+\dots + \frac{\alpha_r}{z-\gamma_r}\right)
\circ W, \ \ \ \ \ D= \left(\beta_0+\frac{\beta_1}{z-\gamma_1}+\dots + \frac{\beta_r}{z-\gamma_r}\right)\circ W.$$ Furthermore, if $\alpha$ is the leading 
coefficient of $A$ and $\beta$ is the leading 
coefficient of $B$ then $\alpha \alpha_i^n=\beta \beta_i^n$, $1\leq i \leq r.$ \qed
\ec
\pr Since $\deg A=\deg B$ it follows from Theorem \ref{polusa} that there exist rational functions $A,B,W$ 
such that $C=\tilde C\circ W,$ $D=\tilde D \circ W$, and all the poles of $\tilde C$ and $\tilde D$ are 
simple (the functions $\tilde C$ and $\tilde D$ obviously have the same set of poles coinciding with the set of poles of the function $A\circ \tilde C=B\circ \tilde D$). Denoting these poles 
by $\gamma_i,$ $1\leq i \leq r,$ we conclude that 
$$\tilde C=\alpha_0+\frac{\alpha_1}{z-\gamma_1}+\dots + \frac{\alpha_r}{z-\gamma_r}
, \ \ \ \ \ \tilde D=\beta_0+\frac{\beta_1}{z-\gamma_1}+\dots + \frac{\beta_r}{z-\gamma_r}$$ 
for some  
$\alpha_i, \beta_i\in \C,$ $0\leq i \leq r$ (in case if $\gamma_i=\infty$ for some $i,$  
$1\leq i \leq r,$ the corresponding terms 
should be changed to $\alpha_iz,$ $\beta_iz$). 

Furthermore, if $\alpha$ (resp. $\beta$) is the leading coefficient of $A$ (resp. $B$) then the leading coefficient of the Laurent expansion of the function $A\circ \tilde C$ (resp. $B\circ \tilde D$) 
near $\gamma_i,$ $1\leq i \leq r,$ 
equals $\alpha \alpha_i^n$ (resp. $\beta \beta_i^n$). Since $A\circ \tilde C=B\circ \tilde D$ this implies 
that for any $i,$ $1\leq i \leq r,$ the equality $\alpha \alpha_i^n=\beta \beta_i^n$ holds. 
\qed

Notice that replacing the rational function $W$ in Corollary \ref{copo} by the function $\mu\circ W,$ where 
$\mu$ is an appropriate automorphism of the sphere, we may assume that
$\gamma_1,\gamma_2,\gamma_3$ are any desired points of the sphere.

Finally, let us mention the following corollary of Theorem \ref{polusa} which generalizes the
corresponding pro\-perty of polynomial decompositions.

\bc \la{zep} Suppose that under assumptions of Theorem \ref{polusa} the function $h$ is a generalized polynomial and $\deg f=\deg g.$ Then $f\circ p\sim g\circ q.$
\ec
\pr Set $x=f^{-1}\{\infty\}$. The conditions of the corollary 
and Theorem \ref{polusa} imply that $\tilde p^{-1}\{x\}$ contains a unique point and the multiplicity of 
this point with respect to $\tilde p$ is one. Therefore $\tilde p$ 
is an automorphism. The same is true for $\tilde q.$ \qed

\section{\la{rita} Ritt classes of rational functions}
As it was mentioned above the first Ritt theorem fails to be true for arbitrary rational functions and it is quite interesting to describe the classes of rational functions for which this theorem remains true. 
In this section we propose an approach to this problem.
This approach is especially useful when a sufficiently complete information about double decompositions of the functions from the corresponding class is available. In particular, our method permits to generalize the first 
Ritt theorem to Laurent polynomials using the classification of their double decompositions.

It is natural to assume that considered classes of rational functions  
possess some property of closeness which is formalized in the following definition. Say that 
a set of rational functions $\f R$ is {\it a closed class} if for any $F\in \f R$
the equality $F=G\circ H$ implies that $G\in \f R$, $H\in \f R.$
For example, rational functions 
for which $$\min_{z\in \C\P^1}\vert F^{-1}\{z\}\vert \leq k,$$ where $k\geq 1$ is a fixed number and
$\vert F^{-1}\{z\} \vert$ denotes the cardinality of the set $F^{-1}\{z\}$, 
form a closed class. We will denote this class by $\f R_k.$ 

Say that two maximal decompositions $\f D,\f E$ of a rational function $F$ are {\it weakly equivalent}
if there exists 
a chain of maximal decompositions $\f F_i$, $1\leq i \leq s,$ of $F$ such that 
$\f F_1=\f D,$ $\f F_s\sim \f E,$ and $\f F_{i+1}$ is obtained from $\f F_i,$ 
$1\leq i \leq s-1,$ by replacing two successive functions $A\circ B$ in $\f F_i$ by new functions $C\circ D$ such that \linebreak $A\circ C=B\circ D.$
It is easy to see that this is indeed an equivalence relation. We will denote this equivalence relation by the symbol $\sim_w$. Say that a closed class of rational functions $\f R$ is {\it a Ritt class} 
if for any $F\in \f R$ any two maximal decompositions of $F$ are weakly equivalent.
Finally, say that a double decomposition \be \la{ebis} H=A\circ C=B\circ D\ee of a rational function $H$ is {\it special} 
if $C,$ $D$ are indecomposable, the pair $A,B$ is reducible, 
and 
there exist no rational functions $\tilde A,$ $\tilde B,$ $U$, $\deg U>1$, such that 
\be \la{kulia} A=U \circ \tilde A, \ \ \ \ B=U \circ \tilde B,\ \ \ \
\tilde A\circ C=\tilde B\circ D.\ee 

For decompositions $$\f A :\ A=A_r\circ A_{r-1}\circ \, ... \, \circ A_1, \ \ \ \ \ \f B:\ B=B_s\circ B_{s-1}\circ \, ... \, \circ B_1$$ of rational functions $A$ and $B$ denote 
by $\f A \circ \f B$ the decomposition
$$ A_r\circ A_{r-1}\circ \, ... \, \circ A_1\circ B_s\circ B_{s-1}\circ \, ... \, \circ B_1$$
of the rational function $A\circ B.$ In case if a rational function $R$ is indecomposable we
will denote the corresponding maximal decomposition by the same letter.

\bt \la{rit1} 
Let $\f R$ be a closed class of rational functions. 
Suppose that for any $P\in \f R$ and any special double decomposition $$P=V\circ V_1=W\circ W_1$$ of $P$ the following condition holds: for any  
maximal decomposition $\f V$ of $V$ and any  
maximal decomposition $\f W$ of $W$ the maximal decompositions 
$\f V\circ V_1$ and $\f W\circ W_1$ of $P$ are weakly equivalent.
Then $\f R$ is a Ritt class.
\et 

\pr For a function $H\in \f R$ denote by $d(H)$ the maximal possible length of a maximal decomposition of $H.$ We use the induction on $d(H)$. 

If $d(H)=1$ then any two maximal decompositions of $H$ are weakly equivalent. So, assume that $d(H)>1$ and let
$$\EuScript H_1:\ H=F_r\circ F_{r-1}\circ \, ... \, \circ F_1, \ \ \ \ \ \f H_2:\ H=G_s\circ G_{s-1}\circ \, ... \, \circ G_1$$ be two maximal decompositions 
of a function $H\in \f R.$ Set 
\be \la{forr} F=F_r\circ F_{r-1}\circ \, ... \, \circ F_2, \ \ \ \ G=G_s\circ G_{s-1}\circ \, ... \, \circ G_2\ee and 
consider the double decomposition
\be \la{decc} H=F\circ F_1=G\circ G_1.\ee If the pair $F,G$ is irreducible then Theorem \ref{rit2} implies that
$\f H_1\sim_w\f H_2$ 
and therefore we must consider only the case when the pair $F,G$ is reducible.

If \eqref{decc} is special then $\f H_1\sim_w\f H_2$ in view of the assumption 
of the theorem. So assume that \eqref{decc} is not special and let $\tilde F,$ $\tilde G,$ $U,$ $\deg U>1,$
be rational functions such that 
$$F=U \circ \tilde F, \ \ \ \ G=U \circ \tilde G, \ \ \ \ \tilde F\circ F_1=\tilde G\circ G_1.$$
Denote by $\hat{\f H}_1, \hat{\f H}_2$ the maximal decompositions \eqref{forr}
of the functions $F$ and $G$ and pick some maximal decompositions 
$$\tilde{\EuScript F}:\ \tilde F=\tilde F_{n}\circ \tilde F_{n-1} \circ \, ... \, \circ \tilde F_1, \ \ \ \ \ \ \  \tilde{\EuScript G}:\  \tilde G=\tilde G_{m}\circ \tilde G_{\tilde m-1} \circ \, ... \, \circ \tilde G_1,
$$
$$\EuScript U:\ U=U_l\circ U_{l-1}\circ \, ... \, \circ U_1$$ of the functions  $\tilde F,$ $\tilde G,$ $U$.

Since $\f R$ is closed $F,G\in \f R$. Furthermore,
$d(F),d(G)<d(H).$ Therefore, the induction assumption implies that 
$$\hat{\f H}_1\sim_w \EuScript U\circ \tilde{\EuScript F}, \ \ \ 
\hat{\f H}_2\sim_w\EuScript U\circ \tilde{\EuScript G}$$
and hence 
\be \la{ebt} {\f H}_1\sim_w \EuScript U\circ \tilde{\EuScript F}\circ F_1, \ \ \ 
{\f H}_2\sim_w\EuScript U\circ \tilde{\EuScript G}\circ G_1.\ee 

Similarly, the function 
$\tilde H=\tilde F\circ F_1=\tilde G\circ G_1$ is contained in 
$\f R$ and $d(\tilde H)<d(H)$. Hence, 
\be \la{ebt1} \tilde{\EuScript F}\circ F_1\sim_w \tilde{\EuScript G}\circ G_1.\ee Now \eqref{ebt} and \eqref{ebt1} imply 
that ${\f H}_1\sim_w {\f H}_2.$ \qed

As an illustration of our approach let us prove the first Ritt theorem.

\bc \la{rit3} The class $\f R_1$ is a Ritt class.
\ec

\pr In view of Theorem \ref{rit1} it is enough to prove that a polynomial $H$ has no special 
double decompositions \eqref{ebis}.
So, assume that the pair $A,$ $B$ in \eqref{ebis} is reducible.
By Corollary \ref{fr2} there exist 
polynomials $A_1,$ $B_1,$ $U,$ $V$ such that 
$$A= A_1\circ W, \ \ \ B=B_1\circ V, \ \ \ 
\deg A_1=\deg B_1>1.$$
Furthermore, Corollary \ref{zep} implies that 
$$A_1\circ (W\circ C)\sim B_1\circ (V\circ D).$$ 
Therefore, equalities \eqref{kulia} hold for  
$$U=A_1,\ \ \ \tilde A=W, \ \ \ \tilde B=\mu\circ V ,$$ and an appropriate $\mu\in \Aut(\C\P^1),$  
and hence \eqref{ebis} is not special. \qed

\section{\la{urrr} Solutions of equations \eqref{2} and \eqref{3}}

In this section we solve equations \eqref{2} and \eqref{3}. 

\bl\la{zc} Let $L_1$, $L_2$ be Laurent polynomials such that the equality
\be \la{raf} L_1\circ z^{d_1}=L_2\circ z^{d_2}\ee
holds for some $d_1,d_2 \geq 1$. Then there exists
a Laurent polynomial $R$ 
such that\be \la{iiuu+} L_1=R\circ z^{D/d_1},\ \ \ \ \ \ L_2=R\circ z^{D/d_2},\ee
where $D=LCM(d_1,d_2).$
\el 

\pr For any subgroup $G$ of $\Aut(\C\P^1)$ the set $k_G$ consisting of rational functions $f$ such that $f\circ \sigma=f$
for all $\sigma\in G$ is a subfield $k_G$ of $\C(z)$. Therefore, by the L\"uroth theorem $k_G$ has the form $k_G=\C(\phi_G(z))$ for some rational function $\phi_G.$

Denote by $F$ the Laurent polynomial defined by 
equality \eqref{raf}.
It follows from \eqref{raf} that $F$ is invariant with respect  
the automorphisms
$\alpha_1\  : \ z\rightarrow  \exp(2\pi i/d_1) z,$
$\alpha_2\  : \ z\rightarrow  \exp(2\pi i/d_2) z.$ Therefore, $F$ is invariant with respect to the automorphism group $G$  generated by 
$\alpha_1,$ $\alpha_2.$ Clearly, $\phi_G=z^{D}$ and hence $F=R\circ z^D$ for some Laurent 
polynomial $R$. Now equalities \eqref{iiuu+} follow from equalities 
$$R\circ z^D=(R\circ z^{D/d_1})\circ z^{d_1} =L_1\circ z^{d_1}, \ \ \
R\circ z^D=(R\circ z^{D/d_2})\circ z^{d_2} =L_2\circ z^{d_2}.\ \ \Box$$   

\vskip 0.2cm

Notice that Lemma \ref{zc} implies that if $A,B,L_1,L_2$ is a solution of equation \eqref{3} then  condition 1) of Theorem \ref{1.1} holds.

\vskip 0.2cm
Set $$D_{ n}=\frac{1}{2}\left(z^n+\frac{1}{z^n}\right).$$ Notice that for any $m\vert n$
$$D_{n}=T_{n/m}\circ D_m= D_{n/m}\circ z^m.$$

\bl \la{erer} Let $F$ be a rational function 
such that $$F(z)=F(1/z)=F(\v z),$$ where $\v$ is a root of unity of order $n\geq 1.$ Then there exists a rational function 
$R$ such that $F=R\circ D_{ n}.$
\el 

\pr Let $G_1$ be a subgroup of $\Aut(\C\P^1)$ generated by 
the automorphism
$\alpha_1\  : \ z\rightarrow \nu z,$ where $\nu =exp(2\pi i/n),$ $G_2$ be a subgroup of $\Aut(\C\P^1)$ generated by 
the automorphism $\alpha_2\ : \ z\rightarrow 
\frac{1}{z}$, and $G_3$ be a subgroup of $\Aut(\C\P^1)$ generated by $\alpha_1$ and $\alpha_2$. It is easy to see that generators of the corresponding invariant fields  
are $\phi_{G_1}=z^{n},$ $\phi_{G_2}=D_1,$ and $\phi_{G_3}=D_{ n}.$
Since $F$
is invariant with respect to $G_1$ and $G_2$ it is  
invariant with respect to $G_3$ 
and therefore $F=R\circ D_{ n}$ for some rational function $R$. \qed

\bl \la{copo1} Let $A,B$ be polynomials of the same degree and $L_1$, $L_2$ be Laurent polynomials such that 
\be \la{tyr} A\circ L_1=B\circ L_2\ee and $A\circ L_1\nsim B\circ L_2.$ 
Then there exist polynomials $w_1,w_2$ of degree one, a root of inity $\nu $, and $a\in \C$
such that 
\be \la{barsu} w_1\circ L_1\circ (az)=D_{ r}, \ \ \ w_2\circ L_2\circ (az)=D_{ r}\circ (\nu  z).\ee 

Furthermore, if a polynomial $A$ and a Laurent polynomial $L$ 
satisfy the equation  \be \la{rere} D_r=A\circ L\ee for some $r\geq 1$ then there exist a polynomial $w$ of degree one, a root of inity $\nu $, and $n\geq 1$ such that \be \la{iui} w\circ L=D_n \circ (\nu  z).\ee
\el

\pr Indeed, it follows from Corollary \ref{copo} that there exist a rational function $W$ and $\alpha_0,\alpha_1, \alpha_2,$ $\beta_0,\gamma\in \C$ such that  
$$L_1=\left(\alpha_0+\alpha_1 z+\frac{\alpha_2}{z}\right)\circ W, \ \ \ \ L_2=\left(\beta_0+\alpha_1\nu _1\gamma z+\frac{\alpha_2\nu _2\gamma}{z}\right)\circ W,$$
for some $r$th roots of unity $\nu _1,\nu _2$.
Furthermore, it follows from $A\circ L_1\nsim B\circ L_2$ that $\alpha_1\alpha_2\neq 0$. Since the function defined by 
equality \eqref{tyr} has two poles this implies that $W=cz^r,$ $c\in \C,$
and without loss of generality we may assume that $c=1.$
The first part of the lemma follows now from the equalities 
$$\alpha_0+\alpha_1 z^r+\frac{\alpha_2}{z^r}=\left(\alpha_0+\frac{2 \alpha_1 z}{a^r} \right)\circ \frac{1}{2}\left(z^r+\frac{1}{z^r}\right)\circ (az),$$ 
$$\beta_0+\alpha_1\nu_1\gamma z^r+\frac{\alpha_2\nu _2\gamma}{z^r}=\left(\beta_0+\frac{2\alpha_1\nu _1 \gamma z}{a^r\nu^r}\right)\circ \frac{1}{2}\left(z^r+\frac{1}{z^r}\right)\circ (\nu az),$$ where $a$ and $\nu$ are numbers satisfying $a^{2r}=\alpha_1/\alpha_2$ and $\nu^{2r}=\nu_1/\nu_2.$ 

Suppose now that equality \eqref{rere} holds. Set $n=\deg L_1$ and consider the equality 
\be \la{les}D_r=T_{r/n}\circ D_n=A\circ L.\ee If the decompositions appeared in \eqref{les} are 
not equivalent then arguing as above 
and taking into account that in this case $a=1$ we conclude 
that \eqref{iui} holds for some root of unity $\nu$. On the other hand, if the decompositions in \eqref{les} are 
equivalent then \eqref{iui} holds for $\nu=1.$ 
\qed

The theorem below provides a description of solutions of equation \eqref{2} and implies that
if $A,L_1,L_2,z^d$ is a solution of \eqref{2} then either condition 1) or condition 4) of Theorem \ref{1.1} holds.

\bt \la{gop} Suppose that polynomials $A,D$ and Laurent polynomials
$L_1,$ $L_2$ (which are not polynomials) satisfy the equation 
\be \la{zx} A\circ L_1=L_2\circ D.\ee
Then there exist polynomials $R,$ $\tilde A,$ $\tilde D,$ 
$W$ and Laurent polynomials $\tilde L_1,$ $\tilde L_2$ such 
that 
\be \la{yyttrr} A=R \circ \tilde A, \ \ \ \ L_2=R\circ \tilde L_2,  \ \ \ \ 
L_1=\tilde L_1 \circ W, \ \ \ \ D= \tilde D\circ W,  \ \ \ \
\tilde A\circ \tilde L_1=\tilde L_2\circ \tilde D \ee
and either 
\be \la{hy} \tilde A\circ \tilde L_1\sim z^{n}\circ z^rL(z^n),\ \ \ \ 
\tilde L_2\circ \tilde D \sim z^{r}L^{n}(z)\circ z^n,\ee
where $L$ is a Laurent polynomial, $r\geq 0,$ $n\geq 1,$ and $\GCD(r,n)=1,$ or 
\be \la{yh} 
\tilde A\circ \tilde L_1\sim  T_{n} \circ D_m,\ \ \ \ 
\tilde L_2\circ \tilde D \sim 
D_m\circ z^n,\ee
where $T_n$ is the $n$th Chebyshev polynomial, $n\geq 1,$ $m \geq 1,$ and $\GCD(m,n)=1.$
\et

\pr Without of loss of generality we may assume that $\C(L_1,D)=\C(z)$. Since 
the function defined by equality \eqref{zx} has two poles, $D=cz^n,$ where $c\in \C,$ 
and we may assume that $c=1.$
Therefore, 
$$A \circ L_1=L_2\circ D=L_2\circ D\circ \v z=
A \circ L_1\circ \v z,$$ where $\v= exp(2\pi i/n).$

If the decompositions $A\circ L_1$ and $A\circ (L_1\circ \v z)$ are equivalent then we have:
\be \la{mu} L_1\circ \v z=\nu \circ L_1,\ee where $\nu \in \Aut(\C\P^1).$ Furthermore, since $\nu$ transforms infinity to infinity, $\nu$ is a linear function and equality \eqref{mu} implies that $\nu^{\circ n}=z$. Therefore, $\nu=\alpha+\omega z$ for some $n$th root of unity $\omega$ and $\alpha\in \C.$
Now the comparison of the coefficients of both parts 
of equality \eqref{mu} 
implies that $L_1$ has the form $$ L_1=\beta+
z^rL(z^n), \ \ \ 0\leq r<n,$$ where $L$ is a Laurent polynomial and $\beta\in \C.$ Clearly, without loss of generality we may assume that $\beta=0$ and this implies that also $\alpha=0.$

It follows from 
$$A \circ L_1=
A \circ L_1\circ \v z=A \circ \omega z\circ L_1$$
that $A\circ \omega z=A$. Since $\omega=\v^r$ and $\GCD(r,n)=1$ in view of the assumption $\C(L_1,D)=\C(z)$,
this implies that  
$ A=R\circ z^{n}$ for some polynomial $R$.
It follows now from the equality
$$
L_2\circ z^{n}=A\circ L_1=R\circ z^{n}\circ z^{r}L(z^{n})=R\circ z^{r}L^{n}(z)
\circ z^{n}$$ 
that
$ L_2=R\circ z^{r}L^{n}(z) .$ Therefore, if the decompositions $A\circ L_1$ and $A\circ (L_1\circ \v z)$ are equivalent then equalities \eqref{yyttrr}, \eqref{hy} hold.

Suppose now that the decompositions $A\circ L_1$ and $A\circ (L_1\circ \v z)$ are not equivalent. Since for any $a\in \C$ we have $z^n\circ (az)=(a^n z)\circ z^n$, it follows from Lemma \ref{copo1} that without loss of generality we may assume that 
$D$ is still equal $z^n$ while 
\be \la{egikk} L_1=D_m=D_1\circ z^m.\ee 
Moreover, $\GCD(m,n)=1$ in view of the assumption $\C(L_1,D)=\C(z)$.
It follows now from \eqref{zx} and \eqref{egikk} and 
Lemmas \ref{zc} and \ref{erer} 
that the Laurent polynomial $L$ defined by
equality \eqref{zx} has the form
$L=R\circ D_{nm},$ where 
$R$ is a polynomial. Therefore, 
$$A\circ D_m=R\circ D_{nm}=
R\circ T_n\circ D_m$$ and hence 
$A=R\circ T_n$. Similarly, 
$$L_2\circ z^n=R\circ D_{nm}=
R\circ D_m\circ z^n$$
and hence $L_2=R\circ D_m.$ \qed

\section{\la{reduc} Reduction of equation \eqref{1} for reducible pairs $A,$ $B$}
In this section we show that the description of solutions of equation \eqref{1} for reducible pairs $A,B$ reduces either to the irreducible case or to the description of
double decompositions of the function $D_n.$ 

\bl \la{posl} Suppose that polynomials $A,B$ satisfy the equation 
\be \la{vtv} A\circ D_{ n}\circ (\mu  z)=B\circ D_{ m},\ee
where $\gcd(n,m)=1$ and $\mu $ is a root of unity.  
Then there exist
a polynomial $R$ and $l\geq 1$ 
such that $\mu ^{2nml}=1$ and $$ A=R\circ \mu ^{nml}T_{lm},\ \ \ \ \ \ B=R\circ T_{ln}.$$ \el

\pr 
Let $F$ be a Laurent polynomial defined by equality \eqref{vtv}. 
It follows from $F=B\circ D_m$ 
that $F\circ (1/z)=F.$ On the other hand, 
$$F\circ (1/z)=A\circ   D_n\circ (\mu /z)=A\circ   \frac{1}{2}\left(\left(\frac{\mu } {z}\right)^n+\left(\frac{z}{\mu }\right)^n\right)=$$
$$=A\circ D_n\circ (z/\mu )=A\circ D_n\circ (\mu z)\circ (z/\mu ^2)=F\circ (z/\mu ^2).$$
Therefore, $F=\tilde F\circ z^d$ for some rational function $\tilde F$ and $d$ equal to the order of $1/\mu ^2$. Since also 
$$D_{ n}\circ (\mu  z)= \frac{1}{2}\left(\mu ^{2n} z+\frac{1}{\mu ^{2n} z}\right)\circ z^{n}, \ \ \ 
D_{ m}=D_1\circ z^{m},$$ 
Lemmas \ref{zc} and \ref{erer} imply that $F=R\circ D_{nml},$
where $R$ is a rational function and $l={\rm lcm}(d,nm)/nm.$

It follows now from 
$$B\circ  D_m=R\circ D_{nml}=R\circ T_{ln}\circ  D_m
$$ that $B=R\circ T_{ln}$. On the other hand, taking into account that 
$\mu ^{nml}=\pm 1,$
we have: 
$$A\circ  D_{ n}=F\circ  (z/\mu )=R\circ D_{nml}\circ (z/\mu )=R\circ \mu ^{nml}D_{nml}=
R\circ \mu ^{nml} T_{lm}\circ  D_n$$ 
and therefore $A=R\circ( \mu ^{nml}T_{lm}).$

\bt \la{red} Suppose that polynomials $A,B$ and Laurent polynomials
$L_1,$ $L_2$ satisfy the equation  
\be \la{eq-}A\circ L_1=B\circ L_2\ee
and the pair $A,B$ is reducible. Then
there exist polynomials $R,$ $\tilde A,$ $\tilde B,$ 
$W$ and Laurent polynomials $\tilde L_1,$ $\tilde L_2$ such 
that 
\be \la{yhy} A=R \circ \tilde A, \ \ \ \ B=R\circ \tilde B,  \ \ \ \ 
L_1=\tilde L_1 \circ W, \ \ \ \ L_2= \tilde L_2\circ W, \ \ \ \ \tilde A\circ \tilde L_1=\tilde B\circ \tilde L_2\ee and either the pair $\tilde A, \tilde B$ is irreducible or
\be \la{yh2+} 
\tilde A\circ \tilde L_1\sim  
- T_{nl}\circ \frac{1}{2}\left(\v z^m+\frac{1}{\v z^m}\right),
\ \ \ \ 
\tilde B\circ \tilde L_2 \sim 
T_{ml} \circ\frac{1}{2}\left(z^n+\frac{1}{z^n}\right),\ee
where $T_{nl},T_{ml}$ are the corresponding Chebyshev polynomials with $n,m\geq 1$, $l>2,$ $\v^{nl}=-1,$
and $\GCD(n,m)=1.$ 
\et

\pr Without loss of generality we may assume that $\C(L_1,L_2)=\C(z)$
and that 
there exist no rational functions $R,$ $\tilde A,$ $\tilde B$ with $\deg R>1$ 
such that the equalities  
\be \la{ffgg}  A=R \circ \tilde A, \ \ \ \ B=R\circ \tilde B,  \ \ \  \tilde A\circ L_1=\tilde B\circ L_2\ee hold.
If the pair $A,B$ is irreducible then there is nothing to prove so assume that it 
is reducible.

By Theorem \ref{p6} and 
Corollary \ref{fr2} there exist polynomials $A_1,$ $B_1,$ $U,$ $V$ such that  
\be \la{tusha} A= A_1\circ U, \ \ \ \ B=B_1\circ V, \ \ \ \deg A_1=\deg B_1>1. \ee
Furthermore,
\be \la{suka} A_1\circ (U\circ L_1)\nsim B_1\circ (V\circ L_2)\ee since 
otherwise \eqref{ffgg} holds for 
$$R=A_1, \ \ \ \tilde A=U,\ \ \ \  \tilde B=\mu\circ V,$$ where $\mu$ is an appropriate automorphism of the sphere.
Therefore, by the first part of Lemma \ref{copo1}, 
we may assume without loss of generality
that  
\be \la{zaii} U\circ L_1 = D_r\circ (\nu z), \ \ \ V\circ L_2=D_r,\ee
where $\nu$ is a root of unity. Applying now the second part of Lemma \ref{copo1} to equalities \eqref{zaii} we see that without loss of generality we may assume that
\be \la{ryryr} L_1=D_m\circ (\mu z),\ \ \ L_2=D_n,\ee
where $\mu$ is a root of unity. Moreover, $\GCD(n,m)=1$ in view of the condition $\C(L_1,L_2)=\C(z)$. In particular, we may assume that $n$ is odd.

It follows from \eqref{ryryr} by Lemma \ref{posl} taking into account the assumption about 
solutions of \eqref{ffgg} 
that there exists a polynomial
$R$ of degree one such that 
$$ A=R\circ (\v^{nl} T_{nl}), \ \ L_1=\frac{1}{2}\left(\v z^m+\frac{1}{\v z^m}\right),\ \ B=R\circ T_{ml},  \ \ L_2=D_n,$$ where $\v=\mu^{m}$ and $l\geq 1$. 
Furthermore, since the pair $A,B$ is reducible it follows from Proposition \ref{p4} that $l>1$. 
Clearly, $\v^{2nl}=1.$ Notice finally that we may assume that $\v^{nl}=-1.$ Indeed, if $\v^{nl}=1$ and $nl$ is 
odd then, taking into account that $T_{nl}\circ (-z)=-T_{nl},$
we may just change $\v$ to $-\v$. On the other hand, if $nl$ is
even then 
$\v^{nl}=1$ contradicts to the assumption about solutions of \eqref{ffgg}.
Indeed, since by the assumption $n$ is odd, if $nl$ is even then $l$ is also even and 
$\v^{nl}=1$ implies that $\mu^{mn(l/2)}=\pm 1$. Hence, $$T_{nl}=T_2\circ (\mu^{mn(l/2)} T_{n(l/2)})$$ and
$$A=(R\circ T_2)\circ ((\mu^{mn(l/2)} T_{n(l/2)})\circ D_m\circ (\mu z))
, \ \ B=(R\circ T_2)\circ T_{m(l/2)}\circ D_n,$$ where
$$(\mu^{mn(l/2)} T_{n(l/2)})\circ D_m\circ (\mu z)=(\mu^{mn(l/2)}D_{mn(l/2)})\circ (\mu z)=D_{mn(l/2)}=T_{m(l/2)}\circ D_n.$$

In order to finish the proof we only must show that the algebraic curve \be \la{tit} T_{ln}(x)+T_{lm}(y)=0,\ee where $GCD(n,m)=1,$ 
is reducible if and only if $l>2$. First observe that if $l$ is divisible by an odd number $f$ then \eqref{tit} is reducible
since $$T_{ln}(x)+T_{lm}(y)=T_f\circ T_{n(l/f)}-T_f\circ( -T_{m(l/f)}).$$ Similarly, if $l$ is divisible by 4 then \eqref{tit} is also reducible since the curve \linebreak $T_4(x)+T_4(y)=0$ is reducible. 

On the other hand, if $l=2$ then \eqref{tit} is irreducible. Indeed, 
otherwise Corollaries \ref{bbv}, \ref{fr2} imply that 
\be \la{gygy} T_{2n}=A_1\circ C, \ \ \ -T_{2m}=B_1\circ D, \ee for some 
polynomials $A_1,B_1,$ $C,$ $D$ such that $\deg A_1=\deg B_1=2$ and the curve 
\be \la{ebkol}  A_1(x)-B_1(y)=0 \ee is reducible. Since $T_{2k}=T_2\circ T_k$ 
it follows from Corollary \ref{zep} that if equalities \eqref{gygy} hold 
then $A_1=T_2\circ \mu_1,$ $B_1=-T_2\circ \mu_2$ for some automorphisms of the sphere.
However, it is easy to see that in this case curve \eqref{ebkol} is not reducible.
Therefore, the condition that equality \eqref{gygy} holds and the condition that curve \eqref{ebkol} is reducible may not be satisfied simultaneously and hence \eqref{tit} is irreducible.
\qed

\section{\la{irreduc} Solutions of equation \eqref{1} for irreducible pairs $A,$ $B$}
In this section we describe solutions of equation \eqref{1} in the case when
the pair $A,B$ is irreducible. We start from a general description of the approach 
to the problem.

First of all, if $A,B$ is an irreducible pair of polynomials then rational functions 
$C,D$ sa\-tisfying equation \eqref{-0} exist if and only if the genus of
curve \eqref{cur} equals zero. Furthermore, it follows from Theorem \ref{p2} that 
if $\tilde C,$ $\tilde D$ is a rational solution of \eqref{-0} 
such that $\deg \tilde C=\deg B,$ $\deg \tilde D=\deg A$ then for 
any other rational solution $C,D$ of \eqref{-0} there exist rational functions $C_1,$ $D_1,$ $W$
such that
$$C=C_1\circ W, \ \ \ D=D_1 \circ W, \ \ \ A\circ C_1\sim A\circ  \tilde C, \ \ \ B\circ D_1\sim B\circ \tilde D.$$
Finally, if $C,D$ are Laurent polynomials then the function $h_1$ from Theorem \ref{p2} should have two poles. On the other hand, it follows from the description of the monodromy of $h_1$, taking into account that 
$A,B$ are polynomials, that the number 
of poles of $h_1$ equals $\GCD(\deg A,\deg B).$

The remarks above imply that in order to describe   
solutions of equation \eqref{1} for irreducible pairs of polynomials $A,B$
we must describe all irreducible pairs of polynomials $A,B$ such that $\GCD(\deg A,\deg B)\leq 2$ and 
the expression for the genus of \eqref{cur} provided 
by formula \eqref{rh0} gives zero. Besides, for each of such pairs we must find a pair of Laurent polynomials $\tilde L_1,$ $\tilde L_2$ satisfying \eqref{1} and such that 
$\deg \tilde L_1=\deg B,$ $\deg \tilde L_2=\deg A.$

The final result is the following statement which supplements (over the field $\C$) 
Theorem 6.1 of the paper
of Bilu and Tichy \cite{bilu}.

\bt \la{irrr} 
Suppose that polynomials $A,B$ and Laurent polynomials $L_1,$ $L_2$ satisfy the equation
$$A\circ L_1=B\circ L_2$$ and the pair $A,B$ is irreducible. 
Then there exist polynomials $\tilde A,$ $\tilde B$, $\mu,$ $\deg \mu=1,$ and 
rational functions $\tilde L_1,$ $\tilde L_2,$ $W$
such that 
$$A=\mu \circ \tilde A, \ \ \ \ B=\mu \circ \tilde B,\ \ \ \    
L_1=\tilde L_1 \circ W, \ \ \ \ L_2=\tilde L_2 \circ W, \ \ \ \ \tilde A\circ \tilde L_1=\tilde B\circ \tilde L_2$$ and, up to a possible replacement of $A$ to $B$ and $L_1$ to $L_2$,
one of the following conditions holds:
$$\tilde A\circ \tilde L_1\sim z^n \circ z^rR(z^n),  \ \ \ \ \ \ \tilde B\circ \tilde L_2\sim  z^rR^n(z) \circ z^n,\leqno 1)$$ 
where $R$ is a polynomial, $r\geq 0,$ $n\geq 1,$ and $\GCD(n,r)=1;$ 
$$\tilde A\circ \tilde L_1\sim T_n \circ T_m, \ \ \ \ \ \ \tilde B\circ \tilde L_2\sim T_m \circ T_n,\leqno 2)$$ 
where $T_n,T_m$ are the corresponding Chebyshev polynomials with $m,n\geq 1,$ and $\GCD(n,m)=1;$ 
$$\tilde A\circ \tilde L_1\sim -T_{2n_1} \circ \frac{1}{2}\left(\v z^{m_1}+\frac{1}{\v z^{m_1}}\right), \ \ \ \ \ \ 
\tilde B\circ \tilde L_2\sim T_{2m_1} \circ  \frac{1}{2}\left(z^{n_1}+\frac{1}{z^{n_1}}\right),\leqno 3)$$ 
where $T_{2n_1},T_{2m_1}$ are the corresponding Chebyshev polynomials with $m_1,n_1\geq 1,$ $\v^{2n_1}=-1,$ and $\GCD(n_1,m_1)=1;$ 
$$\tilde A\circ \tilde L_1\sim  z^2 \circ \frac{z^2-1}{z^2+1}
S(\frac{2z}{z^2+1}),  \ \ \ \ \ \ \tilde B\circ \tilde L_2\sim  (1-z^2)S^2(z)\circ \frac{2z}{z^2+1},\leqno 4)$$
where $S$ is a polynomial;
$$\tilde A\circ \tilde L_1\sim (z^2-1)^3\circ 
\frac{3(3z^4+4z^3-6z^2+4z-1)}{(3z^2-1)^2},  \leqno 5)$$$$
\ \ \ \ \ \ \tilde B\circ \tilde L_2\sim (3z^4-4z^3)\circ \frac{4(9z^6-9z^4+18z^3-15z^2+6z-1)}{(3z^2-1)^3}.$$
\et

The proof of this theorem is given below and consists of the following stages. 
First we rewrite formula for the genus of \eqref{cur} in a more convenient way and
prove several related lemmas. Then we introduce the conception of a special value
and classify the polynomials having such values.
The rest of the proof reduces to the analysis of two cases: the case when one 
of polynomials $A,B$ does not have
special values and the case when both $A,B$ have special values. 

Notice that if at least one of polynomials $A,B$ (say $A$) is of degree 1 then condition 1) holds with 
$\mu =A,$ $R=A^{-1}\circ B,$ $n=1,$ $r=0,$ $W=L_2.$ So, below we always will assume that 
$\deg A,\deg B>1.$ Besides, since one can check by a direct calculation that all the pairs of Laurent polynomials $\tilde L_1,$ $\tilde L_2$ in Theorem \ref{irrr} satisfy the requirements above, we will concentrate on the finding of $A$ and $B$ only.

\subsection{Genus formula and related lemmas}
Let $S=\{z_1,z_2, \dots ,z_s\}$ be any set of complex numbers which contains all {\it finite} 
branch points of a polynomial $A$ of degree $n$.
Then the collection of partitions of the number $n$: $$(a_{1,1},a_{1,2}, ... , a_{1,p_1}),\ \dots \ ,(a_{s,1},a_{s,2}, ... , a_{s,p_{s}}),
$$ where $(a_{i,1},a_{i,2}, ... , a_{i,p_i}),$ $1\leq i \leq s,$ is 
the set of lengths of disjoint cycles in the permutation 
$\alpha_i(A)$, is called the {\it passport} of $A$ and is denoted by $\f P(A)$.
Notice that, since we do not require that any of the points of $S$ is a branch point of $A$, 
some of partitions above may contain units only. We will call such partitions
trivial and will denote by $s(A)$ the number of non-trivial partitions in $\f P(A)$.

Below we will assume that $S$ is a union of all finite branch points of a pair 
of polynomials $A,B$, $\deg A=n$, $\deg B=m,$ and use the notation  
$$(b_{1,1},b_{1,2}, ... , b_{1,q_1}),\ \dots \ ,(b_{s,1},b_{s,2}, ... , b_{s,q_{s}}),
$$ for the passport $\f P(B)$ of $B$.
Clearly, by the Riemann-Hurwitz formula we have:
\be \la{f} 
\sum_{i=1}^{s}{p_i}=(s-1)n+1, \ \ 
\sum_{i=1}^{s}{q_i}=(s-1)m+1.\ee 

For an irreducible pair of polynomials $A,$ $B$ denote by $g(A,B)$ the genus of curve \eqref{cur}. We start from 
giving a convenient version of formula \eqref{rh0} for $g(A,B)$.

\bl \la{11} 
\be \la{rh2} -2g(A,B)=\GCD(m,n)-1+$$$$+\sum_{i=1}^{s}\sum_{j_1=1}^{p_{i}}\left[
a_{i,j_1}(1-q_i)-1
+\sum_{j_2=1}^{q_{i}} \GCD(a_{i,j_1}b_{i,j_2})\right].
\ee
\el

\pr It follows from \eqref{f} that 
$$
\sum_{i=1}^{s}\sum_{j_1=1}^{p_{i}}\left[
a_{i,j_1}(1-q_i)-1\right]= 
\sum_{i=1}^{s}\left[n(1-q_i)-p_i\right]=
ns-n\sum_{i=1}^{s}q_i-\sum_{i=1}^{s}p_i=$$$$
=ns-n((s-1)m+1)-((s-1)n+1)=
-n(s-1)m-1.
$$
Therefore, the right side of formula \eqref{rh2} equals 
$$-n(s-1)m-2+\sum_{i=1}^{s}\sum_{j_1=1}^{p_{i}}
\sum_{j_2=1}^{q_{i}} \GCD(a_{i,j_1}b_{i,j_2})+\GCD(m,n)$$
Now \eqref{rh2} follows from \eqref{rh0} taking into account that $r=s+1.$ \qed

\vskip 0.2cm

Set
$$s_{i,j_1}=a_{i,j_1}(1-q_i)-1+\sum_{j_2=1}^{q_{i}} \GCD(a_{i,j_1}b_{i,j_2}),$$
$1\leq i \leq s,$ $1\leq j_1 \leq p_i.$ 
Using this notation we may rewrite formula \eqref{rh2} in the form
\be \la{rys2} -2g(A,B)=\GCD(m,n)-1+\sum_{i=1}^{s}\sum_{j_1=1}^{p_{i}}s_{i,j_1}.
\ee
Two lemmas below provide upper estimates for $s_{i,j_1},$ $1\leq i \leq s,$ $1\leq j_1 \leq p_i.$

\bl \la{rr-} In the above notation for any fixed pair of indices
$i$, $j_1$, $1\leq i \leq s,$ $1\leq j_1 \leq p_i,$ the following statements hold:
\vskip 0.2cm
\noindent a) If there exist 
at least three numbers $b_{i,l_1},b_{i,l_2},b_{i,l_3},$ $1\leq l_1,l_2,l_3 \leq q_i,$ 
which are not divisible by $a_{i,j_1}$ then 
$s_{i,j_1}\leq -2;$ 
\vskip 0.2cm
\noindent b) If
there exist exactly two numbers $b_{i,l_1},b_{i,l_2},$ $1\leq l_1,l_2 \leq q_i,$ 
which are not divisible by $a_{i,j_1}$ then $s_{i,j_1}\leq -1$  
and the equality attains if and only if 
\be \la{gcdd} \GCD(a_{i,j_1}b_{i,l_1})=\GCD(a_{i,j_1}b_{i,l_2})=a_{i,j_1}/2;\ee

\noindent c) If there exists exactly one number $b_{i,l_1},$ $1\leq l_1\leq q_i,$ which is not divisible by $a_{i,j_1}$
then \be \la{poi} s_{i,j_1}= -1+\GCD(a_{i,j_1}b_{i,l_1}).\ee 
\el

\pr 
If there exist 
at least three numbers $b_{i,l_1},b_{i,l_2},b_{i,l_3},$ $1\leq l_1,l_2,l_3 \leq q_i,$ 
which are not divisible by $a_{i,j_1}$ then 
we have: 
$$s_{i,j_1}=a_{i,j_1}(1-q_i)-1+\sum_{\substack{j_2=1\\j_2\neq l_1,l_2,l_3}}
^{q_{i}}\GCD(a_{i,j_1}b_{i,j_2})
+\sum_{l_1,l_2,l_3}\GCD(a_{i,j_1}b_{i,l_1})\leq $$
$$\leq 
a_{i,j_1}(1-q_i)-1+(q_i-3)a_{i,j_1}+3a_{i,j_1}/2=-a_{i,j_1}/2-1\leq  -2.$$

If
there exist exactly two numbers $b_{i,l_1},b_{i,l_2},$ $1\leq l_1,l_2 \leq q_i,$ 
which are not divisible by $a_{i,j_1}$ then 
we have: 
$$s_{i,j_1}=a_{i,j_1}(1-q_i)-1+\sum_{\substack{j_2=1\\j_2\neq l_1,l_2}}
^{q_{i}}\GCD(a_{i,j_1}b_{i,j_2})
+\sum_{l_1,l_2}\GCD(a_{i,j_1}b_{i,l_1})\leq $$
$$\leq a_{i,j_1}(1-q_i)-1+(q_i-2)a_{i,j_1}+a_{i,j_1}/2+a_{i,j_1}/2=-1,$$
and the equality attains if and only if 
$$\GCD(a_{i,j_1}b_{i,l_1})=\GCD(a_{i,j_1}b_{i,l_2})=a_{i,j_1}/2.$$ 

Finally, if there exists exactly one number $b_{i,l_1}$ which is not divisible by $a_{i,j_1}$
then we have:
$$s_{i,j_1}=a_{i,j_1}(1-q_i)-1+\sum_{\substack{j_2=1\\j_2\neq l_1}}
^{q_{i}}\GCD(a_{i,j_1}b_{i,j_2})
+\GCD(a_{i,j_1}b_{i,l_1})=$$
$$=a_{i,j_1}(1-q_i)-1+(q_i-1)a_{i,j_1}+\GCD(a_{i,j_1}b_{i,l_1})=-1+\GCD(a_{i,j_1}b_{i,l_1}). \ \ \ \Box$$
\vskip 0.2cm

\bc \la{ebanis} Let $B$ be a polynomial of degree $m$ such that the curve $x^n- B(y)=0$ is irreducible and of genus zero. Then  
\vskip 0.2cm
\noindent a)
The equality $\GCD(n,m)=1$ implies that there exists a polynomial $\nu$ of degree 1 such that $B\circ \nu =z^rR^n$ for some polynomial $R$ and $r\geq 1$ such that $\GCD(r,n)=1;$
\vskip 0.2cm
\noindent b) The equality $\GCD(n,m)=2$ implies that $n=2$ and there exists a polynomial $\nu$ of degree 1 such that 
$B\circ \nu=(1-z^2)S^2$ for some polynomial $S$.
\ec
\pr First of all observe that it follows from the irreducibility of
$x^n- B(y)=0$ that among the numbers $b_{1,1},b_{1,2}, ... , b_{1,q_1}$ there exists at least
one number which is not divisible by $n.$

If $\GCD(m,n)=1$ then it follows from formula \eqref{rys2} that $s_{1,1}=0$ and Lemma \ref{rr-} implies that
all the numbers 
$b_{1,1},b_{1,2}, ... , b_{1,q_1}$ but one, say $b_{1,1}$,  
are divisible by $n$ while $\GCD(n,b_{1,1})=1.$
Clearly, this implies that $B\circ \nu =z^rR^n$ for some $\nu,$ $R$, and $r$ 
as above. 

Similarly, if $\GCD(m,n)=2$ then it follows from formula \eqref{rys2} that $s_{1,1}=-1$ and Lemma \ref{rr-} implies that all the numbers 
$b_{1,1},b_{1,2}, ... , b_{1,q_1}$ but two, say 
$b_{1,1},b_{1,2}$,   
are divisible by $n$ while $\GCD(n,b_{1,1})= \GCD(n,b_{1,2})=n/2$. Since this implies that 
$B=z^{n/2}\circ W$ for some polynomial $W$ it follows now from the irreducibility of $x^n- B(y)=0$
that $n=2$ and therefore $B\circ \nu =(1-z^2)S^2$ for some $\nu$ and $S$ 
as above. 
\qed

\bc \la{rr} In the notation of Lemma \ref{rr-} suppose additionally that 
\be \la{gcd} \GCD(b_{i1},b_{i2}, ... , b_{iq_{i}})=1. \ee Then the following statements hold:
\vskip 0.2cm
\noindent a) $s_{i,j_1}\leq 0;$ 
\vskip 0.2cm
\noindent b) $s_{i,j_1}= 0$ if and only if either
$a_{i,j_1}=1$ or  
all the numbers $b_{i,j_2},$ 
$1\leq j_2 \leq q_i,$ except one 
are divisible 
by $a_{i,j_1};$ 
\vskip 0.2cm
\noindent c) 
$s_{i,j_1}= -1$ if and only if $a_{i,j_1}=2$ and all the numbers $b_{i,j_2},$ 
$1\leq j_2 \leq q_i,$ but two 
are even.
\ec

\pr If $a_{i,j_1}=1$ then $s_{i,j_1}=0$ so assume that $a_{i,j_1}>1.$ 
The assumption \eqref{gcd} implies that among the numbers $b_{1,1},b_{1,2}, ... , b_{1,q_1}$ there exists at least
one number which is not divisible by $n.$
If there exists exactly one number $b_{i,l_1}$ which is not divisible by $a_{i,j_1}$
then in view of \eqref{gcd} necessarily $\GCD(a_{i,j_1}b_{i,l_1})=1$ and hence $s_{i,j_1}=0$
by formula \eqref{poi}.
If
there exist exactly two numbers $b_{i,l_1},b_{i,l_2},$ $1\leq l_1,l_2 \leq q_i,$ 
which are not divisible by $a_{i,j_1}$ then it follows from Lemma \ref{rr-} that 
$s_{i,j_1}\leq -1$ where the equality attains if and only if \eqref{gcdd} holds. On the other hand, if \eqref{gcdd} holds then 
necessarily $a_{i,j_1}=2$ since otherwise we obtain a contradiction with \eqref{gcd}. 
Finally, if there exist 
at least three numbers $b_{i,l_1},b_{i,l_2},b_{i,l_3},$ $1\leq l_1,l_2,l_3 \leq q_i,$ 
which are not divisible by $a_{i,j_1}$ then $s_{i,j_1}\leq -2$ by Lemma \ref{rr-}. \qed

\subsection{Polynomials with special values} In the notation above say 
that $z_i,$\linebreak $1\leq i \leq s,$ is a special value of $B$ if 
$$\GCD(b_{i,1},b_{i,2}, ... , b_{i,q_{i}})>1.$$ 
It is easy to see that a polynomial $P$ has a special 
value if and only if there exists $c\in \C$ such that 
$P-c=z^d\circ R$ for some polynomial $R.$

Say that $z_i,$ $1\leq i \leq s,$ is a $1$-special value (resp. a $2$-special value) 
of $B$ if all the numbers 
$$b_{i,1},b_{i,2}, ... , b_{i,q_{i}}$$ 
but one (resp. two) are divisible by some number $d>1.$

\bp \la{ch} Let $B$ be a polynomial. Then the following statements hold:
\vskip 0.2cm
\noindent a) $B$ may not have two special values, or one special value and one 1-special value, or
three $1$-special values;
\vskip 0.2cm
\noindent b) If $B$ has two $1$-special values then $s(B)=2,$  
$\f P(B)=\{(1,2,\, \dots\, 2),\ (1,2,\, \dots\, 2)\};$ 
\vskip 0.2cm
\noindent c) If $B$ has one 1-special value and one 2-special value then $s(B)=2$ and either 
$\f P(B)=\{(1,1,2),\ (1,3)\}$
or $\f P(B)=\{(1,2,2),\ (1,1,3)\}.$
\ep

\pr Set $m=\deg B.$
Suppose first that $B$ has at least two 1-special 
values. To be definite assume that these values are $z_1,z_2$ and that
all $(b_{1,1}, \dots ,b_{1,q_1})$ but
$b_{1,1}$ are divisible by the number $d_1$ and 
all $(b_{2,1}, \dots ,b_{2,q_2})$ but
$b_{2,1}$ are divisible by the number $d_2.$ 
Then 
\be \la{zi} q_{1}\leq 1+\frac{m-b_{1,1}}{d_1}, \ \ \ \ \  
q_{2}\leq 1+\frac{m-b_{2,1}}{d_2},\ee
where the equalities attain if only if
$b_{1,j}=d_1$ for $1<j\leq q_1$
and $b_{2,j}=d_2$ for $1<j \leq q_2.$
Furthermore, clearly \be \la{1q} \sum_{i=1}^{s}{q_i}\leq 
q_1+q_2
+(s-2)m,\ee where the equality attains if and only if 
$(b_{i,1}, \dots ,b_{i,q_i})=(1,1,\, \dots\, 1)$
for any $i>2.$ Finally, for $i=1,2$ we have:
\be \la{2q} q_i\leq 1+\frac{m-b_{i,1}}{d_i}\leq 1 +\frac{m-1}{2}     \ee
and hence 
\be \la{2qq} q_1+q_2\leq 1 +m,\ee
where the equality attains only if 
$d_1=2$, $d_2=2,$ $b_{1,1}=1,$ $b_{2,1}=1.$
Now \eqref{1q} and \eqref{2qq} imply that
\be\la{kll} \sum_{i=1}^{s}{q_i}\leq (s-1)m+1.\ee
Since however in view of \eqref{f} in this inequality 
should attain equality we conclude that in all intermediate inequalities should attain equalities.
This implies that $s(B)=2$ and
$$(b_{1,1}, \dots ,b_{1,q_1})=(1,2,\, \dots\, 2),\ \ \ \ (b_{2,1}, \dots ,b_{2,q_1})=(1,2,\, \dots\, 2).$$ 
In particular, we see that $B$ may not have three $1$-special values.
 
In order to prove the first part of the proposition it is enough to observe that if for at least one index 1 or 2, say 1,
the corresponding point is special then 
$$q_1\leq \frac{m}{d_1}\leq \frac{m}{2}.$$ Since this inequality is stronger than \eqref{2q} repeating 
the argument above we obtain an inequality in \eqref{kll} in contradiction with \eqref{f}.

Finally, assume that $z_1$ is a 1-special value while $z_2$
is a 2-special value. We will suppose that all $(b_{1,1}, \dots ,b_{1,q_1})$ but
$b_{1,1}$ are divisible by the number $d_1$ and all
$(b_{2,1}, \dots ,b_{2,q_2})$ but
$b_{2,1},b_{2,2}$ are divisible by the number $d_2.$ 

If $m$ is odd then $d_2\neq 2.$ Hence, in this case  
$d_2\geq 3,$ 
$$ q_{1}\leq 1+\frac{m-b_{1,1}}{d_1}\leq 1+\frac{m-1}{2}, \ \ \ \ \  
q_{2}\leq 2+\frac{m-b_{2,1}-b_{2,2}}{d_2}\leq 2+\frac{m-2}{3},$$ and, therefore,
$$q_1+q_2\leq \frac{11}{6}+\frac{5m}{6}.$$ If $m>5$ then 
$$q_1+q_2\leq \frac{11}{6}+\frac{5m}{6}<m+1.$$ 
Since combined with \eqref{1q} the last inequality leads to a contradiction with \eqref{f}
we conclude that $m\leq 5.$ It follows now from $d_2\geq 3$ that necessarily $m=5$ and 
$(b_{2,1}, \dots ,b_{2,q_2})=(1,1,3).$ Finally, since $z_1$ is a 1-special value
of $B$ we necessarily have 
$(b_{1,1}, \dots ,b_{1,q_1})=(1,2,2).$   
 
Similarly, if $m$ is even then $d_1\geq 3$ and we have:
$$ q_{1}\leq 1+\frac{m-b_{1,1}}{d_1}\leq 1+\frac{m-1}{3}, \ \ \ \ \  
q_{2}\leq 2+\frac{m-b_{2,1}-b_{2,2}}{d_2}\leq 2+\frac{m-2}{2},$$ and
$$q_1+q_2\leq \frac{5}{3}+\frac{5m}{6}\ .$$
If $m>4$ then 
$$\frac{5}{3}+\frac{5m}{6}<m+1$$ and as above we obtain a contradiction with \eqref{f}. 
On the other hand, if $m\leq 4$ then $d_1\geq 3$ implies that necessarily $m=4$ and $(b_{1,1}, \dots ,b_{1,q_1})=(1,3)$. Finally, clearly $(b_{2,1}, \dots ,b_{2,q_2})=(1,1,2).$ \qed

\subsection{Proof of Theorem \ref{irrr}. Part 1} First of all notice that
if at least one of polynomials $A,B$ 
has a unique finite branch point or equivalently is of the form $\mu \circ z^d \circ \nu$ for some $d\geq 1$ and polynomials $\mu,\nu$ of degree one
then 
it follows from Corollary \ref{ebanis} that either condition 1) or condition 4) of Theorem \ref{irrr} holds. So, below we always will assume that both polynomials
$A$, $B$ have at least two finite branch points.
In this subsection we prove Theorem \ref{irrr} under the assumption that at least 
one of polynomials $A,B$ does not have special values. Without loss of
generality we may assume that this polynomial is $B.$ In other words, we may assume that for any $i,$ $1\leq i \leq s,$ equality \eqref{gcd} holds.

\vskip 0.2cm
\noindent
{\it Case 1.} Suppose first that $\GCD(n,m)=1$.  
In this case by formula \eqref{rys2} the condition $g(A,B)=0$ 
is equivalent to the condition 
\be \la{rh+} \sum_{i=1}^{s}\sum_{j_1=1}^{p_{i}}s_{i,j_1}=0.
\ee 
In view of Corollary \ref{rr},a this is possible if and only $s_{i,j_1}=0,$ $1\leq i \leq s,$ 
$1\leq j_1 \leq p_i.$

Since $A$
has at least two finite branch points, Corollary \ref{rr},\,a and Corollary \ref{rr},\,b, taking into account that $B$ may not have more than two 1-special values by Proposition \ref{ch},\,a, 
imply that $A$ has exactly two branch points. Furthermore, it follows from 
Proposition \ref{ch},\,b that $\f P(B)$ equals
\be \la{chebysh1}\{(1,2,2, ... ,2),(1,2,2, ... ,2)\}. \ee
Now Corollary \ref{rr},\,b implies that
\be \la{iu} a_{1,j_1}\leq 2, \ \ \ a_{2,j_2}\leq 2, \ \ \ 1\leq j_1 \leq p_1,\ \ \ 1\leq j_2 \leq p_2.\ee
Since $$p_1+p_2=(s-1)n+1=n+1$$ and
$$\sum_{j_1=1}^{p_1}a_{1,j_1}+\sum_{j_1=1}^{p_2}a_{2,j_1}=2n$$ it follows 
from \eqref{iu} that among $a_{1,j_1},$ $a_{2,j_2}$, $1\leq j_1 \leq p_1,$ $1\leq j_2 \leq p_2,$ there are exactly two units and therefore $\f P(A)$ 
equals either \eqref{chebysh1}  
or 
\be \la{chebysh2}
\{(1,1,2, ... ,2),(2,2,2, ... ,2)\}.\ee 

Recall that for any polynomial $P$ such that $\f P(P)$ equals \eqref{chebysh1} or \eqref{chebysh2} there exist polynomials $\mu,$ $\nu$ of degree 1 such that $\mu\circ P \circ \nu=T_n$ for some $n\geq 1$.
A possible way to establish it is to observe that it follows from $T_n(\cos z)=\cos nz$ that $T_n$ satisfies the differential equation 
\be \la{difur} n^2(y^2-1)=(y^{\prime})^2(z^2-1), \ \ \ y(1)=1.\ee On the other hand, it is easy to see that if $\f P(P)$ equals \eqref{chebysh1} or \eqref{chebysh2} and $\deg P=n$ then $P$ satisfies the equation  
$$n^2(y-A)(y-B)=(y^{\prime})^2(z-a)(z-b),$$ for some $A,B,a,b \in \C$ 
with $y(b)=A$ or $B$. Therefore for appropriate 
polynomials $\mu,$ $\nu$ of degree 1 the polynomial $\mu\circ P \circ \nu$ 
satisfies the equation \eqref{difur} and hence $\mu\circ P \circ \nu=T_n$ by the 
uniqueness theorem for solutions of differential equations.

Since $\f P(B)$ equals \eqref{chebysh1} and 
$\f P(A)$ equals either \eqref{chebysh1} or
\eqref{chebysh2} 
the above characterization of Chebyshev polynomials implies now
that if $\GCD(n,m)=1$ then condition 2) holds.

\vskip 0.2cm
\noindent
{\it Case 2.} If $\GCD(n,m)=2$ then the condition $g(A,B)=0$ 
is equivalent to the condition that one number from $s_{i,j_1},$ $1\leq i \leq s,$ 
$1\leq j_1 \leq p_i,$ equals -1 while others equal 0.

Since $A$ has at least two branch points, 
Corollary \ref{rr},\,b and Corollary \ref{rr},\,c, taking into account that if $B$ has two 1-special values then $B$ does not have 2-special values by
Proposition \ref{ch},\,b,
imply that $A$ has two branch points and $B$ has one 1-special value and one 2-special value.
Therefore, by Proposition \ref{ch},\,c, $\f P(B)$ equals either 
\be \la{pizden1} \f \{(1,2,2),(1,1,3)\}\ee or
\be \la{pizden2} \f \{(1,3),(1,1,2)\}\ee 
Furthermore, since the assumption 
$\GCD(n,m)=2$ implies that $\deg B$ is even we conclude that 
$\f P(B)$ necessarily equals \eqref{pizden2}. 
It follows now from Corollary \ref{rr},\,b and Corollary \ref{rr},\,c
that for any $j_1,$ $1\leq j_1 \leq p_1,$ 
the number $a_{1,j_1}$ equals 1 or 3 and the partition $(a_{2,1},a_{2,2},\,\dots\,a_{2,p_2})$ 
contains one element
equal 2 and others equal 1. 

Denote by $\alpha$ (resp. by $\beta$) the number of appearances of 1 (resp. of 3) in the first partition of 
$\f P(A)$ and by $\gamma$ the number of appearances of 1 in the second partition of 
$\f P(A)$.
We have:
$$\alpha+3\beta=n, \ \ \   2+\gamma=n,$$
and, by \eqref{f}  
\be \la{opop} \alpha+\beta+\gamma=n.\ee The second and the third of the equations above imply that $\alpha+\beta=2$. Hence the partition $(a_{1,1},a_{2,2},\,\dots\,a_{1,p_1})$
is either $(1,3)$ or $(3,3)$ and $\gamma=n-2$ implies  
that either \be \la{pizdos} \f P(A)=\f P(B)=\{(1,3),(1,1,2)\}\ee or \be\la{lys2}\f P(A)=\{(3,3),(2,1,1,1,1)\}.\ee 

Observe now that for any polynomial $R$ for which $\f P(R)$ equals \eqref{pizden2} 
the derivative of $R$ has the form $R^{\prime}=c(z-a)^2(z-b),$ $a,b,c\in \C.$ 
Therefore, there exist polynomials $\mu, \nu$ of degree 1 such that $$\mu\circ R \circ \nu=\int 12z^2(z-1)\d z=3z^4-4z^3.$$ 
Since $A$ and $B$ have the same set of critical values this implies in particular that if
\eqref{pizdos} holds then $A=B\circ \lambda$ for some polynomial $\lambda$ of degree 1 in contradiction with the irreducibility of the curve $A(x)-B(y)=0$. On the other hand, it is easy to see that if equality \eqref{lys2} holds then there exist polynomials $\mu, \nu_1,\nu_2$ of degree 1 such that $$\mu\circ A \circ \nu_1=(z^2-1)^3, \ \ \ \mu\circ B \circ \nu_2=3z^4-4z^3.$$  
Therefore, if $\GCD(n,m)=2$ then condition 5) holds.

\subsection{Proof of Theorem \ref{irrr}. Part 2}
Suppose now that both polynomials $A$ and $B$ have special values. Then 
by Proposition \ref{ch},\,b each of them has a unique special value. 
The special values of $A$ and $B$
either coincide or are different. If they are different then 
\be \la{op} A=(z^{d_1}+\beta_1)\circ \hat A, \ \ \  \ B=(z^{d_2}+\beta_2)\circ \hat B,\ee
for some $\beta_1,\beta_2\in \C,$  
$\beta_1\neq \beta_2,$ and $d_1,d_2>1.$ 
Since the pair
$A$, $B$ is irreducible and $g(A, B)=0$ the pair 
$A_0=z^{d_1}+\beta_1,$ $B_0=z^{d_2}+\beta_2$ is also irreducible and 
\be \la{ert} g(A_0,B_0)=0.\ee
   
Formula \eqref{rh2} implies that   
\be\la{rod} -2g(A_0, B_0)=d_1+d_2-d_1d_2+\GCD(d_1,d_2)-2.\ee
If $\GCD(d_1,d_2)=1$ then \eqref{ert} is equivalent to the equality
$(d_1-1)(1-d_2)=0$ which is impossible.
On the other hand, if $\GCD(d_1,d_2)=2$ then \eqref{ert} is equivalent
to the equality $(d_1-1)(1-d_2)=-1$
which holds if and only if 
$d_1=d_2=2.$ 

Repeatedly using Theorem \ref{p6} and 
Corollary \ref{fr2} we can find  
polynomials $P,$ $Q,$ $U,$ $V$ such that  
$$\hat A= P \circ U, \ \ \ \ \hat B=Q\circ V,\ \ \ \ \ \deg P=\deg Q,$$ and the pair $U,$ $V$ is irreducible.
Setting \be \la{surok} A_1=A_0\circ P, \ \ \  B_1=B_0 \circ Q\ee we see that equality
\eqref{tusha} holds. Furthermore, equivalence \eqref{suka} is impossible since otherwise  
$A_1=B_1\circ \mu$ for some polynomial $\mu$ of degree 1 and it follows from
Corollary \ref{zep} and equalities \eqref{surok} that $A_0=B_0\circ \nu$ for some polynomial $\nu$ of degree 1
in contradiction with the irreducibility of 
the pair $A_0,B_0$.
Now using the same reasoning as in the proof of Theorem \ref{red} and taking into account that the pair 
$A_0,B_0$ is irreducible
we conclude that condition 3) holds. 

\vskip 0.2cm
In the case when the special values of $A$ and $B$
coincide we can assume without loss of generality that 
\be \la{po} A=z^{d_1}\circ U, \ \ \ \ \ \ B=z^{d_2}\circ V,\ee where
$$d_1=\GCD(a_{1,1}, a_{1,2}, \dots, a_{1,p_1})>1, \ \ \ \ d_2=\GCD(b_{1,1}, b_{1,2}, \dots, b_{1,q_1})>1,$$
and \be \la{bg} \GCD(d_1,d_2)=1\ee in view of the irreducibility of the pair $A$ and $B$. 
Notice that, since $A$ and $B$ have 
at least two critical values, the inequalities $p_1\geq 2,$ $q_1\geq 2$ hold. 
Finally, without loss of generality we may assume that $m=\deg B$ is greater than $n=\deg A.$
We will consider the cases $\GCD(d_1,m)=2$ and $\GCD(d_1,m)=1$ se\-parately and will show that in both cases there exist 
no irreducible pairs $A,B$ with $g(A,B)=0.$

\vskip 0.2cm
\noindent
{\it Case 1.} Suppose first that $\GCD(d_1,m)=2.$ Then necessarily $\GCD(n,m)=2$ and,  
since \be \la{vb} x^{d_1}-B(y)=0\ee is an irreducible curve of genus zero, Lemma \ref{ebanis}
implies that $d_1=2$ and all the numbers 
$b_{1,1},b_{1,2}, ... , b_{1,q_1}$ but two, say $b_{1,q_1-1},b_{1,q_1},$ are even while 
$b_{1,q_1-1},b_{1,q_1}$ are odd.
Since by the assumption 
each $a_{1,j_1},$ $1\leq j_1\leq p_1,$ 
is divisible by $d_1=2$, this implies in particular that for each $j_1,$ $1\leq j_1\leq p_1,$ 
$$\GCD(a_{1,j_1}b_{1,q_1-1})\leq a_{1,j_1}/2,\ \ \ \GCD(a_{1,j_1}b_{1,q_1})\leq a_{1,j_1}/2.$$ 

Returning now to polynomials $A,$ $B$ we conclude 
that for each $j_1,$ $1\leq j_1\leq p_1,$ $$ s_{1,j_1}=a_{1,j_1}(1-q_1)-1+\sum_{j_2=1}^{q_{1}-2} 
\GCD(a_{1,j_1}b_{1,j_2})
+\GCD(a_{1,j_1}b_{1,q_1-1})+\GCD(a_{1,j_1}
b_{1,q_1})\leq $$ $$\leq a_{1,j_1}(1-q_1)-1+a_{1,j_1}(q_1-2)+\GCD(a_{1,j_1}b_{1,q_1-1})+\GCD(a_{1,j_1}
b_{1,q_1})\leq $$ 

$$ \leq -a_{1,j_1}-1+a_{1,j_1}/2+a_{1,j_1}/2\leq -1.$$  
Since $p_1\geq 2$ and by Corollary \ref{rr},a for any $i,$ $1 < i \leq s,$ and 
$j_,$ $1\leq j_1 \leq p_i,$ the inequality 
$s_{i,j_1}\leq 0$ holds 
it follows now from formula \eqref{rys2} 
that $g(A,B)<0.$ 

\vskip 0.2cm
\noindent
{\it Case 2}. Similarly, if $\GCD(d_1,m)=1$ then 
Lemma \ref{ebanis} applied to curve 
\eqref{vb} implies that 
each $b_{1,j_1},$ $1\leq j_1\leq q_1,$ except one, say $b_{1,q_1},$ is divisible by $d_1$
while $\GCD(b_{1,q_1},d_1)=1$ and returning to 
$A,B$ and taking into account that 
each $a_{1,j_1},$ $1\leq j_1\leq p_1,$ 
is divisible by $d_1$
we obtain that  
$$s_{1,j_1}= a_{1,j_1}(1-q_1)-1
+\sum_{j_2=1}^{q_{1}-1} 
\GCD(a_{1,j_1}b_{1,j_2})+\GCD(a_{1,j_1}
b_{1,q_1})
\leq $$ \be \la{kot} \leq -1+\GCD(a_{1,j_1}b_{1,q_1})\leq 
-1+a_{1,j_1}/d_1.\ee Hence, 
\be \la{kota} \sum_{j_1=1}^{p_1}s_{1,j_1}\leq -p_1+n/d_1.\ee

Furthermore, since each $b_{1,j_2},$ $1\leq j_2\leq q_1,$ except one is divisible by $d_1$, each 
$b_{1,j_2},$ $1\leq j_2\leq q_1,$ is divisible by $d_2$, and equality \eqref{bg} holds
we have:
$$(q_1-1)d_1d_2+d_2\leq m$$ and therefore 
$$ q_1\leq 1+ m/d_1d_2-1/d_1.$$
Since by \eqref{f} the inequality \be \la{xsw} q_1+q_i\geq m+1\ee holds for any $i,$ $2\leq i \leq s,$ this implies that \be \la{mn} q_i\geq m-m/d_1d_2+1/d_1
.\ee 
Denote by $\gamma_i,$ $2\leq i \leq s,$
the number of units among the numbers $b_{i,j_2},$ $1\leq j_2\leq q_i.$
Since the number of non-units is $\leq m/2$ the inequality $\gamma_i\geq q_i-m/2$ holds and therefore   
\eqref{mn} implies that
\be \la{fg} \gamma_i\geq m/2-m/d_1d_2+1/d_1.\ee 

For any $i,j_1,$ $2\leq i \leq s,$ $1\leq j_1\leq p_i,$ we have:
\be \la{kj} s_{i,j_1} \leq a_{i,j_1}(1-q_i)-1+ a_{i,j_1}(q_i-\gamma_i)+\gamma_i=(1-\gamma_i)(a_{i,j_1}-1).\ee Since 
this implies that
\be \la{krot} \sum_{j_1=1}^{p_i}s_{i,j_1}=(1-\gamma_i)\sum_{j_1=1}^{p_i}(a_{i,j_1}-1)
 \leq  (1-\gamma_i)
(n-p_i)\ee it follows now from \eqref{fg} that 
$$ \sum_{j_1=1}^{p_i}s_{i,j_1}
\leq (1-1/d_1+m(1/d_1d_2-1/2))(n-p_i).$$ 
Therefore, using \eqref{f} we obtain that
\be \la{kota2}\sum_{i=2}^s \sum_{j_1=1}^{p_i}s_{i,j_1} 
\leq(1-1/d_1+m(1/d_1d_2-1/2))(p_1-1).\ee

Set $$S=\sum_{i=1}^s\sum_{j_1=1}^{p_1}s_{i,j_1}.$$ Since $\GCD(n,m)=1$ or 2 it follows from formula \eqref{rys2} that in order 
to finish the proof 
it is enough to show that $S<-1.$  

Since $p_1\geq 2$ it follows from \eqref{kota}, \eqref{kota2} that 
$$S\leq -p_1+n/d_1+
(1-1/d_1+m(1/d_1d_2-1/2))(p_1-1)=$$
$$ =-1+n/d_1-\frac{p_1-1}{d_1}+m(1/d_1d_2-1/2)(p_1-1)<$$ \be \la{kota3}
<-1+n/d_1+m(1/d_1d_2-1/2)(p_1-1).
\ee

If $p_1\geq 3$ then it follows from \eqref{kota3}, taking into account the assumption $m\geq n$ and the inequality $1/d_1d_2- 1/2<0$,
that 
$$S< -1+n(1/d_1+2/d_1d_2-1).$$  
Since $1/d_1+2/d_1d_2-1\leq 0$ for any $d_1,d_2\geq 2,$ this implies that $S<-1.$

If $p_1=2$ then \eqref{kota3} implies that
$$S<-1+n(1/d_1+1/d_1d_2-1/2).$$ 
Since $1/d_1+1/d_1d_2-1/2\leq 0$ whenever $d_1>2$
we obtain that $S<-1$ also if $p_1=2$ but 
$d_1>2$.
 
Finally, if $p_1=2$, $d_1=2$ but $m\geq (3/2)n$ then  it follows from 
equality \eqref{kota3} that
$$S< -1
+n(3/4d_2-1/4).$$ 
Since $d_1=2$ implies $d_2\geq 3$ in view of \eqref{bg}, we conclude again that $S<-1.$

Therefore, the only case when the proof of the inequality $S<-1$ is still not finished is the one when 
$p_1=2$, $d_1=2$, and $n\leq m< (3/2)n.$ In this case apply the reasoning above to $A$ and $B$ switched
keeping the same notation. In other words, assume that
$q_1=2$, $d_2=2$, and \be \la{ws} 2n/3< m \leq n.\ee
Then by \eqref{xsw} we have    
$q_i\geq m-1$ for any $i,$ $2\leq i \leq s.$
Therefore, corresponding partitions of $m$ are either trivial or have the form $(1,1,\dots,1,2)$ and hence  
\be \la{ebun} \gamma_i\geq m-2,\ \ \ 2\leq i \leq s.\ee
It follows now from \eqref{krot}, \eqref{ebun}, \eqref{f},
and \eqref{ws} that
\be \la{vvgg} \sum_{i=2}^s \sum_{j_1=1}^{p_i}s_{i,j_1}  \leq (3-m)(p_1-1)< (3-2n/3)(p_1-1)\leq 3-2n/3.\ee
Since $d_2=2$ implies $d_1\geq 3$ in view of \eqref{bg}, 
it follows now from \eqref{vvgg} and $p_1\geq 2$ that 
$$S < -p_1+n/d_1+ 3-2n/3\leq 1 + n/d_1-2n/3\leq  1-n/3.$$ If $n\geq 6$ then this inequality implies that $S<-1.$
On the other hand, the inequality $n\leq 5$ is impossible since otherwise equalities 
\eqref{po}, \eqref{ws}, and $p_1\geq 2$ imply that $d_1=d_2=2$ in contradiction with \eqref{bg}.

In order to finish the proof of Theorem \ref{irrr} it is enough to notice that
for any choice of $\tilde A,$ $\tilde B$ in conditions 1)-5)
the curve
\be \la{svin} \tilde A(x)-\tilde B(y)=0\ee is indeed irreducible. For cases 1) and 2) this is a corollary of 
Proposition \ref{p4}. For case 3) this was proved in the end of the proof of Theorem \ref{red}. 
In case 4) corresponding curve 
\eqref{svin} is irreducible since otherwise Corollary \ref{bbv}
would imply that there exists a polynomial $T$ such that $\tilde B=z^2\circ T$ in contradiction with $\tilde B=(1-z^2)S^2$. Finally, since $\tilde B$ in 5) is indecomposable 
it follows from Corollary \ref{bbv} taking into account Corollary \ref{fr2} 
that corresponding curve \eqref{svin} is irreducible.

\section{\la{last} Proof of Theorem \ref{1.1}} 
Since the description of double decompositions of 
functions from $\f R_2$ reduces to the corresponding problem for Laurent polynomial
and any double decomposition of  
a Laurent polynomial is equivalent to \eqref{1}, \eqref{2}, or \eqref{3} the first part of Theorem 1.1 follows from Theorem \ref{gop}, Theorem \ref{red}, Theorem \ref{irrr} and Lemma \ref{zc}. 
The proof of the second part is given below.

\bt \la{rittt} The class $\f R_2$ is a Ritt class.
\et

\pr We will use Theorem \ref{rit1} and the first part of Theorem \ref{1.1}. First observe that 
the first part of Theorem \ref{1.1} implies that if 
$A\circ C=B\circ D$ is a double decomposition of a function from $\f R_2$ such that $C,D$ are indecomposable 
and there exist no rational functions $\tilde A,$ $\tilde B,$ $U$, $\deg U>1,$ such that \eqref{kulia} holds
then there exist automorphisms of the sphere $\mu,$ $W$ and 
rational functions $\tilde A,$ $\tilde B,$ $\tilde C,$ $\tilde D $
such that one of  
conclusions of Theorem \ref{1.1} holds.
Moreover, it was shown above that in cases 1)-3) and 6) the pair $\tilde A,$ $\tilde B$ is irreducible. 

Observe now that in case 4) the pair $\tilde A,$ $\tilde B$ is also irreducible. Indeed, since $\GCD(n,m)=1$ it follows from the construction given in Subsection \ref{fsys} that for the pair $f=\tilde A,$ $g=\tilde B$ the permutation $\delta_i,$ $1\leq i \leq r,$ corresponding to the loop around the infinity contains two cycles. Therefore, if the pair $\tilde A,$ $\tilde B$ is reducible then $o(f,g)=2$ and both functions $h_1,h_2$ from Theorem \ref{p2} 
have a unique pole. On the other hand, the last statement contradicts to the fact that 
$h_1=\tilde B\circ v_1,$  $h_2=\tilde B\circ v_2$ for some rational function $v_1,$ $v_2$ since $\tilde B$ has two poles.  
 
Finally, as it was observed in the end of the proof of Theorem \ref{red}, in case 5) the pair $\tilde A,$ $\tilde B$ is reducible whenever $l>2$. Since in this case $\tilde C$ and $\tilde D$ are decomposable unless $n=1,$ $m=1$,
it follows now from Theorem \ref{rit1} that in order to prove the proposition it is enough to check that
for any choice of maximal decompositions
$$-T_{l}=u_d\circ u_{d-1} \circ \dots \circ u_1,  \ \ \ T_{l}=v_l\circ v_{l-1} \circ \dots \circ v_1,$$  
the decompositions 
\be \la{decomp} u_d\circ u_{d-1} \circ \dots \circ u_1\circ \frac{1}{2}\left(\v z+\frac{1}{\v z}\right),\ \ \ \ v_l\circ v_{l-1} \circ \dots \circ v_1\circ \frac{1}{2}\left(z+\frac{1}{z}\right),\ee
where $\v^l=-1,$  
are weakly equivalent. 

Since $T_l=T_d\circ T_{l/d}$ for any $d\vert l,$ it follows from Corollary \ref{zep}  
that any maximal decomposition of $T_{l}$ is equivalent to
$T_l=T_{d_1}\circ T_{d_2} \circ \dots \circ T_{d_s},$
where $d_1, d_2 \, \dots \, d_s$ are prime divisors of $l$ such that $d_1 d_2 ... d_s=l.$ 
Taking into account that for $d\geq 1$
$$T_d\circ \frac{1}{2}\left(z+\frac{1}{z}\right)=\frac{1}{2}\left(z+\frac{1}{z}\right)\circ z^d ,$$
this implies easily that 
both decompositions \eqref{decomp} are weakly equivalent to some decomposition of the form  
$$\frac{1}{2}\left(z+\frac{1}{z}\right) \circ z^{d_1}\circ z^{d_2} \circ \dots \circ z^{d_s}. \ \ \ \Box$$

\bibliographystyle{amsplain}

\begin{thebibliography}{10}

\bibitem {az} R. Avanzi, U. Zannier, 
\textit{The equation $f(X) = f(Y)$ in rational functions $X = X(t)$, $Y = Y(t)$},
Compos. Math. 139, No. 3, 263-295 (2003)

\bibitem {bilu1} Y. Bilu, \textit{Quadratic factors of $f(x)-g(y)$}, Acta Arith. 90, No. 4, 341-355 (1999).

\bibitem {bilu} Y. Bilu, R. Tichy, 
\textit{The Diophantine equation $f(x) = g(y)$},
Acta Arith. 95, No.3, 261-288 (2000).
\bibitem {y} M. Briskin, N. Roytvarf, Y. Yomdin,
{\it Center conditions at infinity for Abel differential equation}, Annals
of Math., to appear.

\bibitem {cg} J. M. Couveignes, L. Granboulan, \textit{Dessins from a geometric point of view},
in The Grothendieck theory of dessins d'enfants, 79-113. Cambridge University Press, 1994.

\bibitem {en} H. Engstrom, \textit{Polynomial substitutions,} Amer. J. Math. 63, 249-255 (1941).

\bibitem {er} A. Eremenko, \textit{Some functional equations connected with the iteration of rational functions}, Leningrad Math. J. 1,  no. 4, 905--919 (1990). 

\bibitem {gs} J. Gutierrez, D. Sevilla,
\textit{Building counterexamples to generalizations for rational functions of Ritt's decomposition theorem},
J. Algebra 303, No. 2, 655-667 (2006).

\bibitem {f1} M. Fried, 
\textit{On a theorem of Ritt and related diophantine problems}, 
J. Reine Angew. Math. 264, 40-55 (1973).

\bibitem {f2} M. Fried,
\textit{Fields of definition of function fields and a problem in the reducibility of polynomials in two variables}, Ill. J. Math. 17, 128-146 (1973).

\bibitem {f3} M. Fried, \textit{
Arithmetical properties of function fields. II. The generalized Schur problem,} Acta Arith. 25, 225-258 (1974)

\bibitem {klein} F. Klein, 
\textit{Lectures on the icosahedron and the solution of equations of the fifth degree},
New York: Dover Publications, (1956).

\bibitem {kur} A. Kurosch, \textit{
The theory of groups. Vol. I,} New York: Chelsea Publishing Company, (1955).

\bibitem {zv} N. Magot, A. Zvonkin, \textit{
Belyi functions for Archimedean solids}, 
Discrete Math. 217, No.1-3, 249-271 (2000).

\bibitem {mas} W. Massey, 
\textit{Algebraic topology: An introduction},
Graduate Texts in Mathematics, Vol. 56, (1981).

\bibitem{mir} R. Miranda, \textit{
Algebraic curves and Riemann surfaces}, 
Graduate Studies in Mathematics. 5. Providence, RI: AMS, American Mathematical Society, (1995).

\bibitem{p} M. Muzychuk, F. Pakovich, \textit{Solution of the polynomial moment problem},
preprint, arXiv:0710.4085.

\bibitem {mp2} M. Muzychuk, F. Pakovich, {\it On maximal decompositions of rational functions}, preprint, arXiv:0712.3869.


\bibitem{p2} F. Pakovich, \textit{On the functional equation $F(A)=G(B)$, where A,B are polynomial and $F,G$ are continuous functions}, Math. Proc. Cam. Phil. Soc.,  143, No. 2, 469-472 (2007).

\bibitem{p1} F. Pakovich, \textit{On polynomials sharing preimages of compact sets and related questions}, Geom. Funct. Anal, 18, No. 1, 163-183 (2008).


\bibitem{pppp} F. Pakovich, \textit{The algebraic curve $P(x)-Q(y)=0$ and 
functional equations,} Complex Var. Elliptic Equ., to appear, arXiv:0804.0736v2. 


\bibitem {ar} F. Pakovich, {\it On analogues of Ritt theorems for rational functions with at most two poles}, Russ. Math. Surv. 63, No. 2, 181-182 (2008).

\bibitem {pr1} F. Pakovich, {\it On the equation $P(f)=Q(g),$ where $P,Q$ are polynomials and $f,g$ are entire functions,} preprint, arXiv:0804.0739.

\bibitem {r1} J Ritt, 
\textit{Prime and composite polynomials},
American M. S. Trans. 23, 51-66 (1922).

\bibitem {r2} J Ritt,
\textit{Equivalent rational substitutions},
American M. S. Trans. 26, 221-229 (1924).

\bibitem {r3} J Ritt,
\textit{Permutable rational functions,}
American M. S. Trans. 25, 399-448 (1923).


\bibitem {sch} A. Schinzel, \textit{Polynomials with special regard to
reducibility}, Encyclopedia of Mathematics and Its Applications
\textbf{77}, Cambridge University Press, 2000.

\bibitem {tor} P. Tortrat, \textit{Sur la composition des polyn\^omes},
Colloq. Math. 55, No.2, 329-353 (1988).

\bibitem {za} U. Zannier, 
\textit{Ritt's second theorem in arbitrary characteristic,}
J. Reine Angew. Math. 445, 175-203 (1993)

\end{thebibliography}

\end{document}